\documentclass[final]{siamart0516}

\usepackage{graphicx}
\usepackage{amsfonts}
\usepackage{amsmath}
\usepackage{amssymb}
\usepackage{amsopn}
\usepackage{algorithm}
\usepackage{algorithmic}
\usepackage{array}
\usepackage{multicol}
\usepackage{booktabs}
\usepackage{bm}
\usepackage{placeins}
\usepackage{capt-of}
\usepackage{cite}

\hypersetup{
pdfinfo={
Title={Generalized Accelerated Composite Gradient Method},
Author={Mihai I. Florea and Sergiy A. Vorobyov},
Subject={A Generalized Accelerated Composite Gradient Method: Uniting Nesterov's Fast Gradient Method and FISTA},
Keywords={estimate sequence, Fast Gradient Method, FISTA, monotone, line-search, composite objective, large-scale optimization}
}
}

\DeclareMathOperator*{\argmin}{arg\,min}
\DeclareMathOperator{\sgn}{sgn}

\newcommand{\cmdaab}[1]{\hat{#1}_{k + 1}}
\newcommand{\cmdaac}[1]{\hat{\bm{#1}}_{k + 1}}
\newcommand{\cmdaad}[1]{^{\rm{#1}}}
\newcommand{\cmdaae}[1]{\in \{1, ... , #1\}}
\newcommand{\cmdaaf}[1]{\textbf{#1.}}

\newcommand{\cmdaah}[2]{line~\ref{#2} of Algorithm~\ref{#1}}
\newcommand{\cmdaai}[2]{line~\ref{#2} of Pattern~\ref{#1}}
\newcommand*\samethanks[1][\value{footnote}]{\footnotemark[#1]}

\newcommand{\cmdaaj}{
\begin{table}[t!]
\centering
\small
\caption{Oracle functions of the five test problems}
\label{labelaab}
\begin{tabular}{lllll} \toprule
& $f(\bm{x})$ & $\Psi(\bm{x})$ & $\nabla f(\bm{x})$ & $prox_{\tau \Psi}(\bm{x})$ \\ \midrule
LASSO & $\frac{1}{2}\|\bm{A} \bm{x} - \bm{b} \|_2^2$ & $\lambda_1 \|\bm{x}\|_1$ & $\bm{A}^T(\bm{A}\bm{x} - \bm{b})$ & $\mathcal{T}_{\tau \lambda_1} (\bm{x})$ \\
NNLS & $\frac{1}{2}\|\bm{A} \bm{x} - \bm{b} \|_2^2$ & $\sigma_{\mathbb{R}_{+}^n}(\bm{x})$ & $\bm{A}^T(\bm{A}\bm{x} - \bm{b})$ & $(\bm{x})_{+}$ \\
L1LR & $\mathcal{I}(\bm{A}\bm{x}) - \bm{y}^T \bm{A} \bm{x}$ & $\lambda_1 \|\bm{x}\|_1$ & $\bm{A}^T(\mathcal{L}(\bm{A}\bm{x}) - \bm{y})$ & $\mathcal{T}_{\tau \lambda_1} (\bm{x})$ \\
RR & $\frac{1}{2}\|\bm{A} \bm{x} - \bm{b} \|_2^2$ & $\frac{\lambda_2}{2} \|\bm{x}\|_2^2$ & $\bm{A}^T(\bm{A}\bm{x} - \bm{b})$ & $\frac{1}{1 + \tau \lambda_2} \bm{x}$ \\
EN & $\frac{1}{2}\|\bm{A} \bm{x} - \bm{b} \|_2^2$ & $\lambda_1 \|\bm{x}\|_1 +
\frac{\lambda_2}{2} \|\bm{x}\|_2^2$ & $\bm{A}^T(\bm{A}\bm{x} - \bm{b})$ & $\frac{1}{1 + \tau \lambda_2} \mathcal{T}_{\tau \lambda_1} (\bm{x})$ \\ \bottomrule
\end{tabular}
\end{table}
}

\newcommand{\cmdaak}{
\begin{figure}[t!] \centering \footnotesize
\begin{minipage}[t]{0.48\linewidth} \centering
\includegraphics[width=\textwidth]{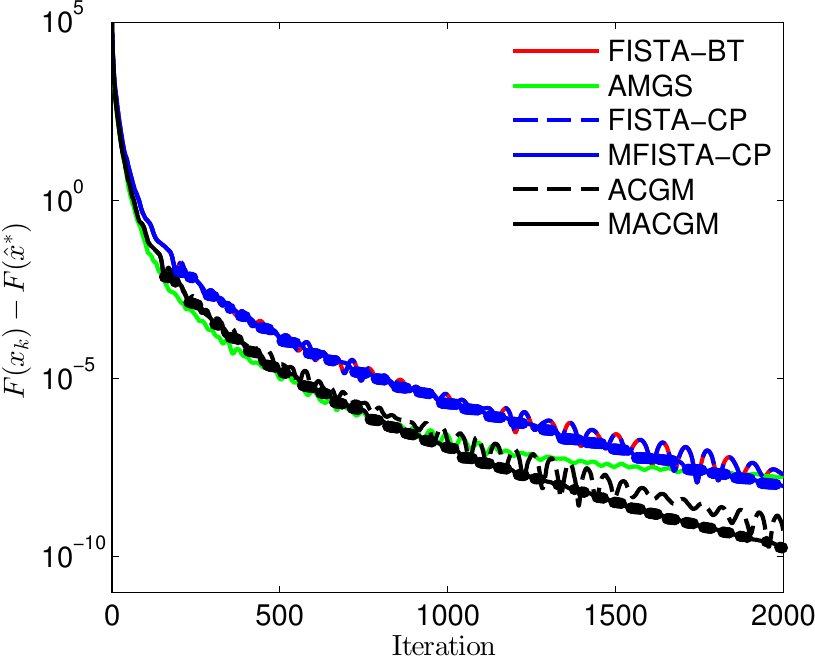}
(a) Iteration convergence rates on LASSO
\end{minipage}
\hspace{2mm}
\begin{minipage}[t]{0.48\linewidth} \centering
\includegraphics[width=\textwidth]{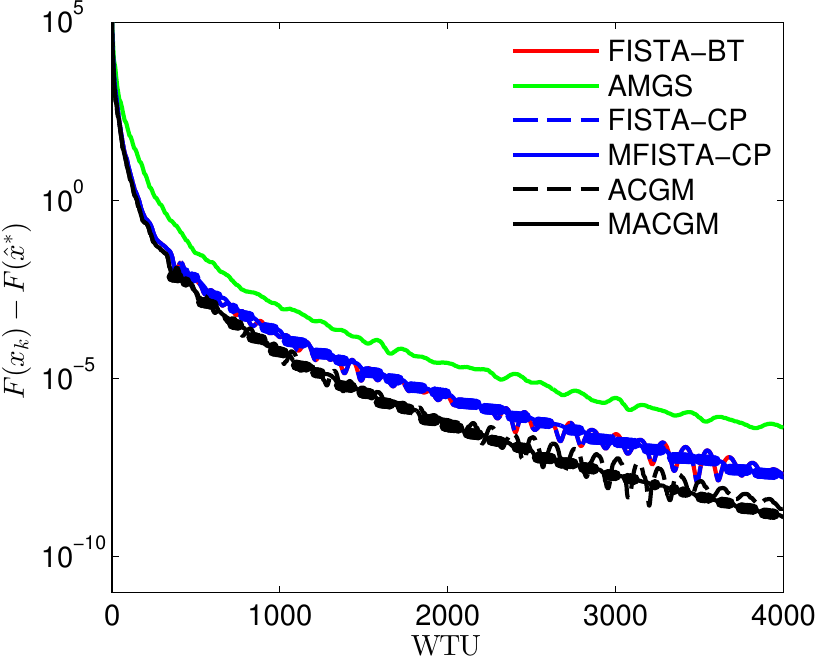}
(b) Computational convergence rates on LASSO
\end{minipage}

\vspace{6mm}

\begin{minipage}[t]{0.48\linewidth} \centering
\includegraphics[width=\textwidth]{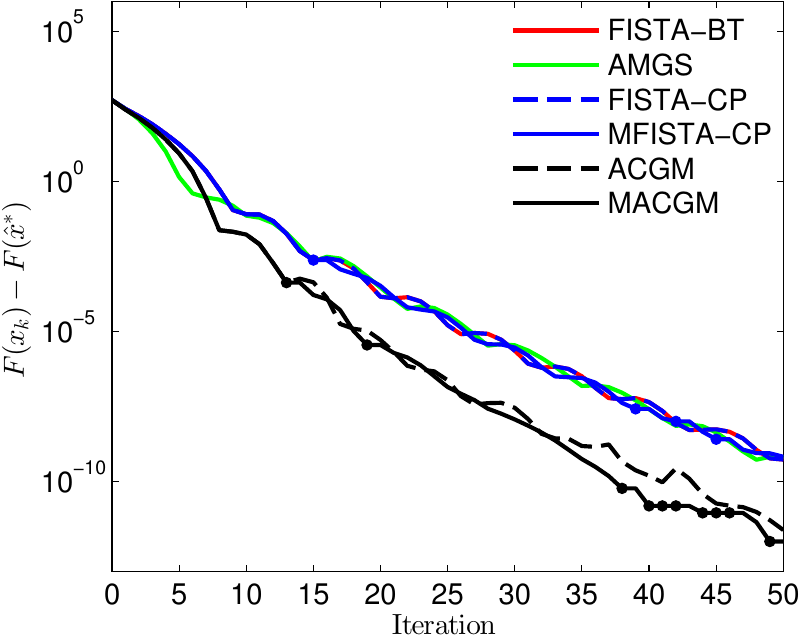}
(c) Iteration convergence rates on NNLS
\end{minipage}
\hspace{2mm}
\begin{minipage}[t]{0.48\linewidth} \centering
\includegraphics[width=\textwidth]{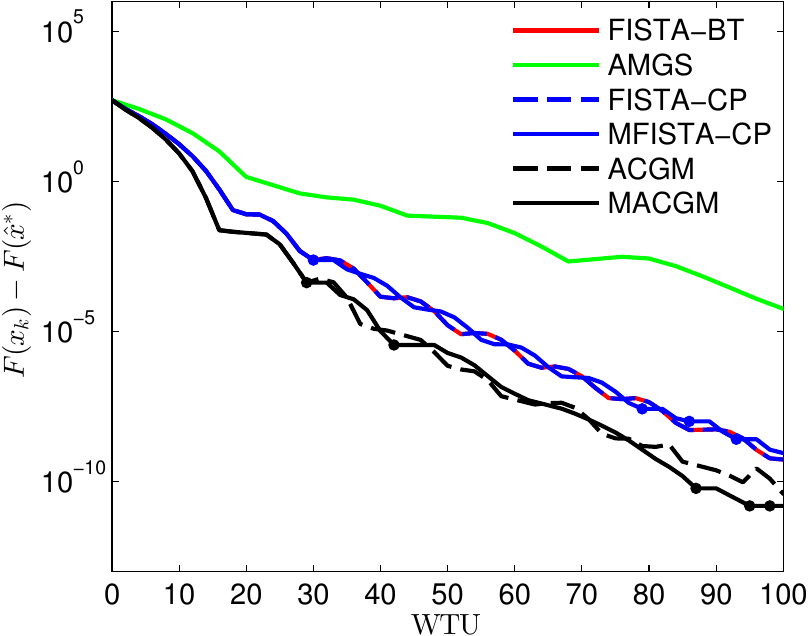}
(d) Computational convergence rates on NNLS
\end{minipage}

\vspace{6mm}

\begin{minipage}[t]{0.48\linewidth} \centering
\includegraphics[width=\textwidth]{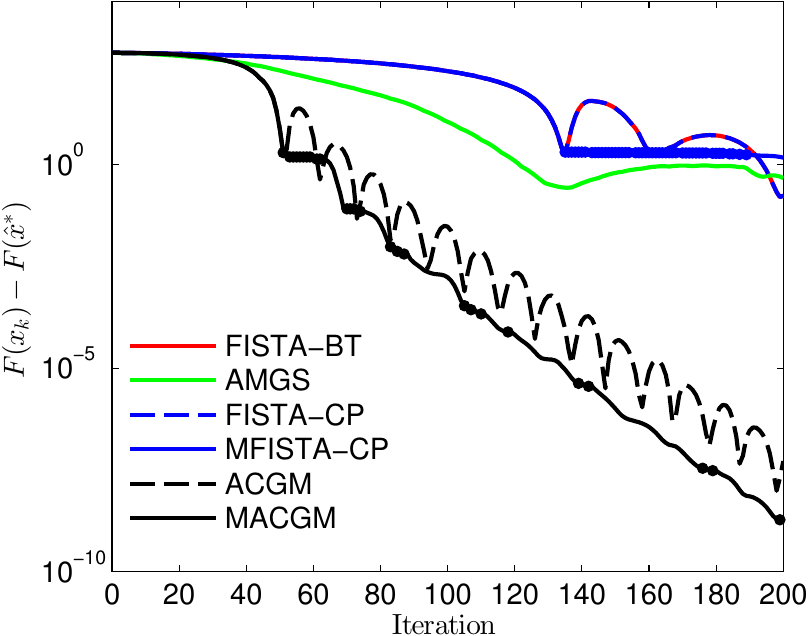}
(e) Iteration convergence rates on L1LR
\end{minipage}
\hspace{2mm}
\begin{minipage}[t]{0.48\linewidth} \centering
\includegraphics[width=\textwidth]{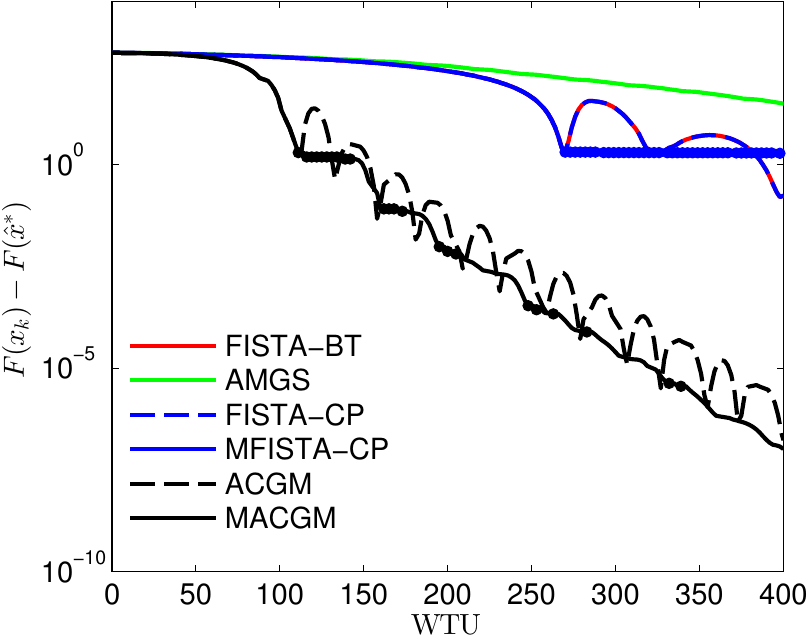}
(f) Computational convergence rates on L1LR
\end{minipage}

\vspace{5mm}

\caption{Convergence results of FISTA with backtracking (FISTA-BT), AMGS, FISTA-CP, monotone FISTA-CP (MFISTA-CP), non-monotone ACGM and monotone ACGM (MACGM) on the LASSO, NNLS, and L1LR non-strongly convex problems. Dots mark iterations \emph{preceding} overshoots. At these iterations, the convergence behavior changes.}
\label{labelaac}
\end{figure}
}

\newcommand{\cmdaal}{
\begin{figure}[t!] \centering \footnotesize
\begin{minipage}[t]{0.49\linewidth} \centering
\includegraphics[width=\textwidth]{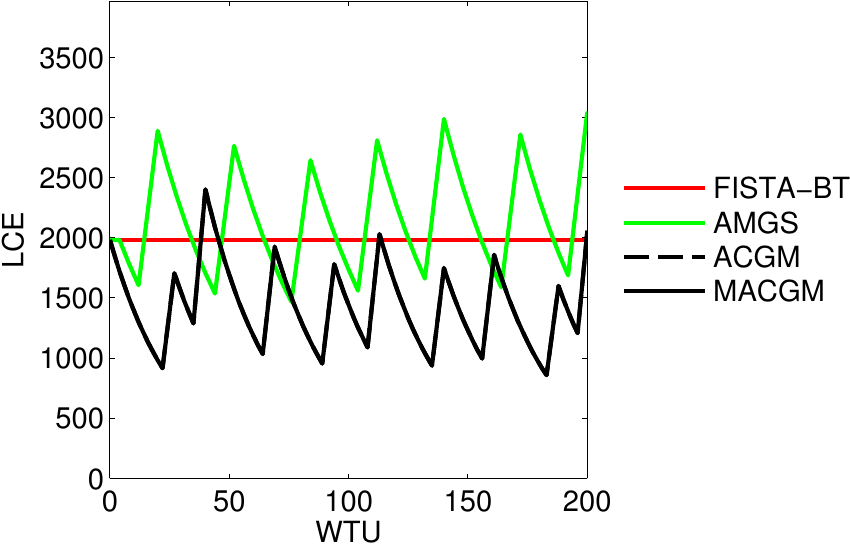}
(a) LASSO
\end{minipage}
\hspace{2mm}
\begin{minipage}[t]{0.47\linewidth} \centering
\includegraphics[width=\textwidth]{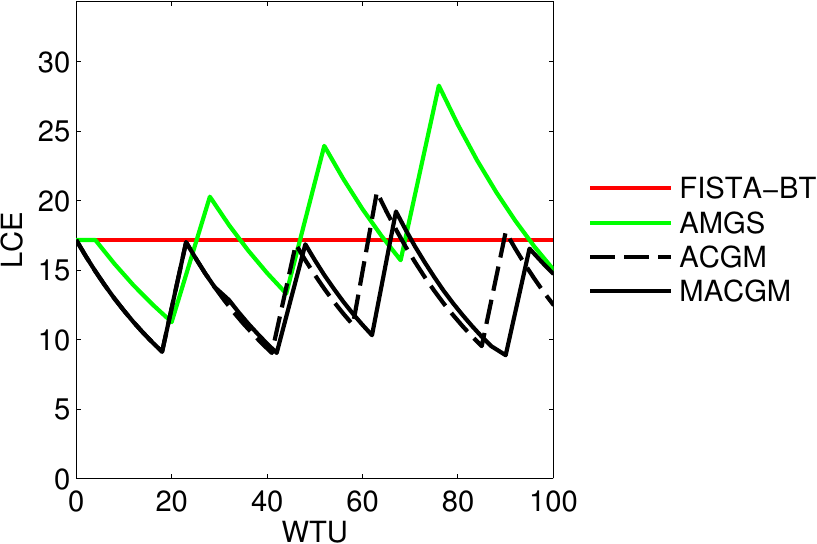}
(b) NNLS
\end{minipage}

\vspace{3mm}

\begin{minipage}[t]{0.48\linewidth} \centering
\includegraphics[width=\textwidth]{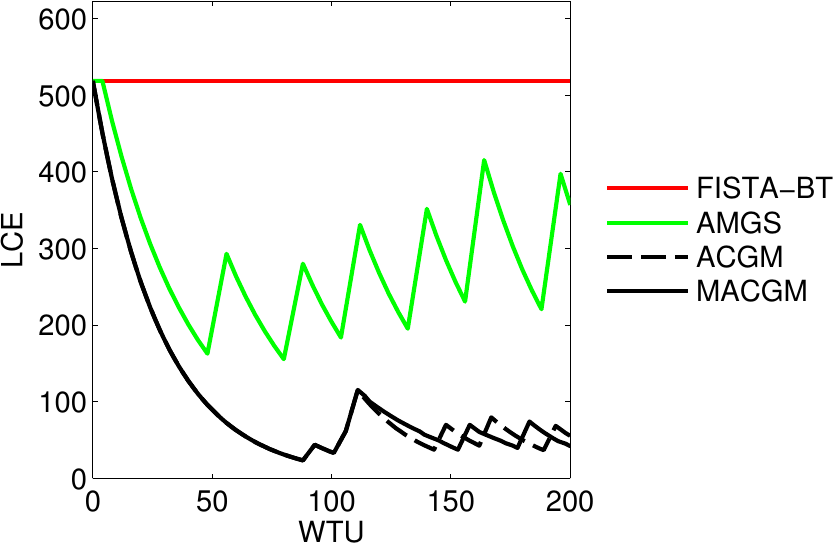}
(c) L1LR
\end{minipage}
\caption{Line-search method LCE variation on LASSO, NNLS, and L1LR}
\label{labelaad}

\vspace{3mm}

\captionof{table}{Average LCEs of line-search methods on LASSO, NNLS, and L1LR}
\label{labelaae}
\small
\begin{tabular}{@{}ccccccc@{}} \toprule
Problem & $L_f$ & Iterations & FISTA-BT & AMGS & ACGM & MACGM \\ \midrule
LASSO & 1981.98 & 2000 & 1981.98 & 2202.66 & 1385.85 & 1303.70 \\
NNLS & 17.17 & 50 & 17.17 & 19.86 & 14.35 & 13.54 \\
L1LR & 518.79 & 200 & 518.79 & 246.56 & 80.76 & 79.12 \\ \bottomrule
\end{tabular}
\end{figure}
}

\newcommand{\cmdaam}{
\begin{figure}[t!] \centering \footnotesize
\begin{minipage}[t]{0.48\linewidth} \centering
\includegraphics[width=\textwidth]{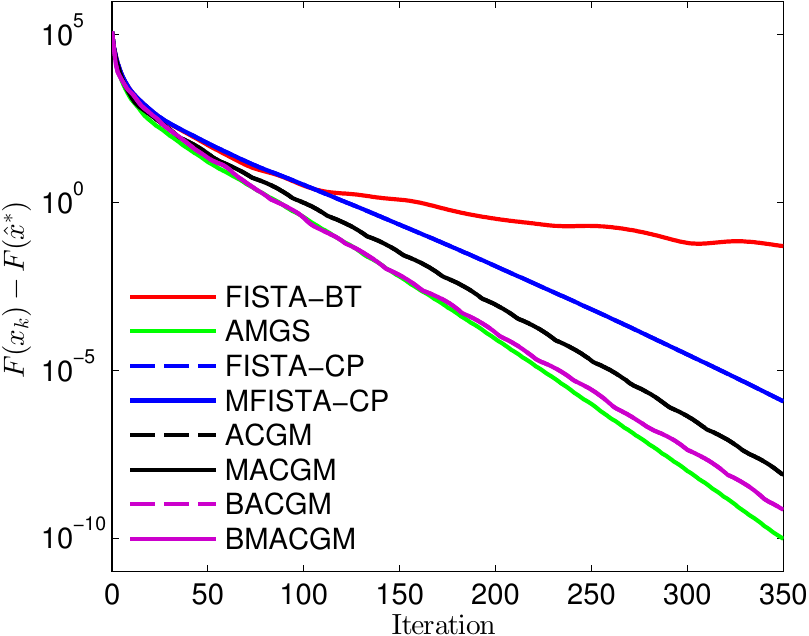}
(a) Iteration convergence rates on RR
\end{minipage}
\hspace{2mm}
\begin{minipage}[t]{0.48\linewidth} \centering
\includegraphics[width=\textwidth]{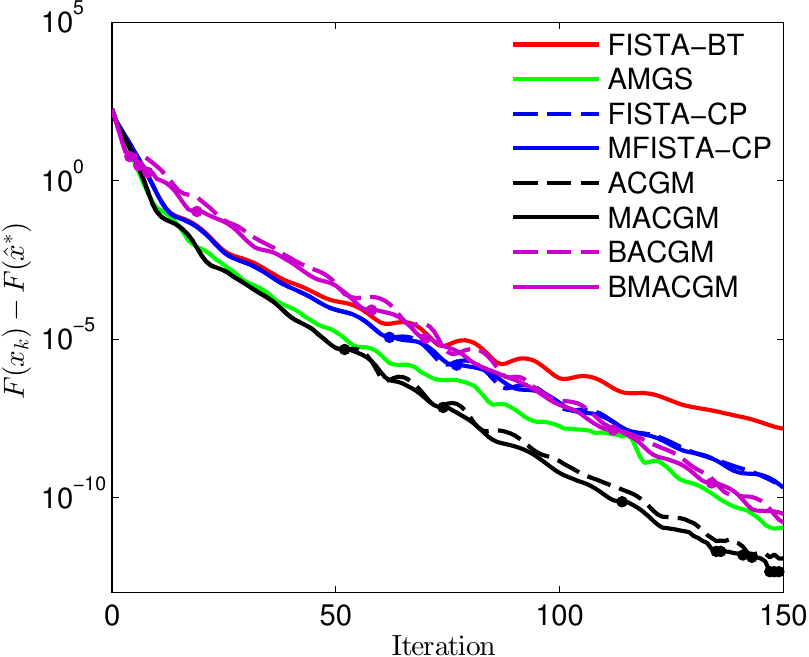}
(b) Iteration convergence rates on EN
\end{minipage}

\vspace{6mm}

\begin{minipage}[t]{0.48\linewidth} \centering
\includegraphics[width=\textwidth]{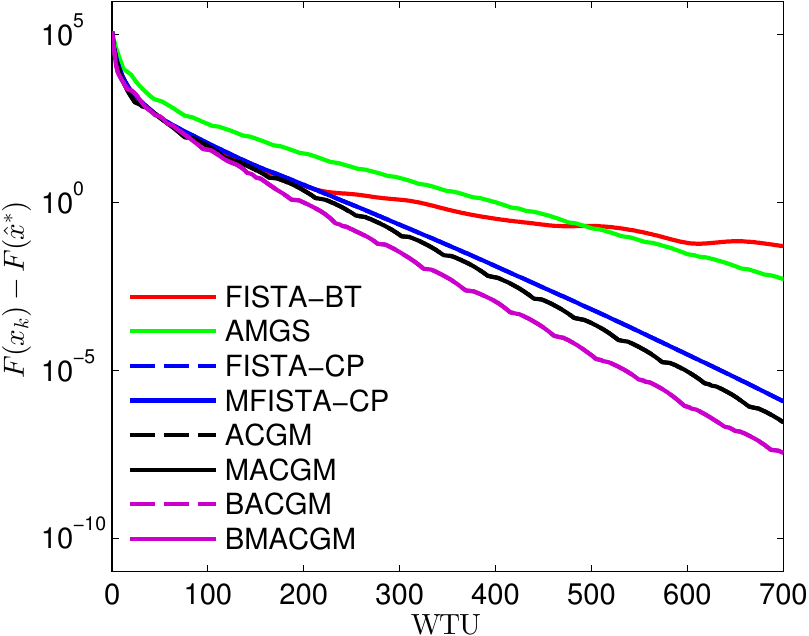}
(c) Computational convergence rates on RR
\end{minipage}
\hspace{2mm}
\begin{minipage}[t]{0.48\linewidth} \centering
\includegraphics[width=\textwidth]{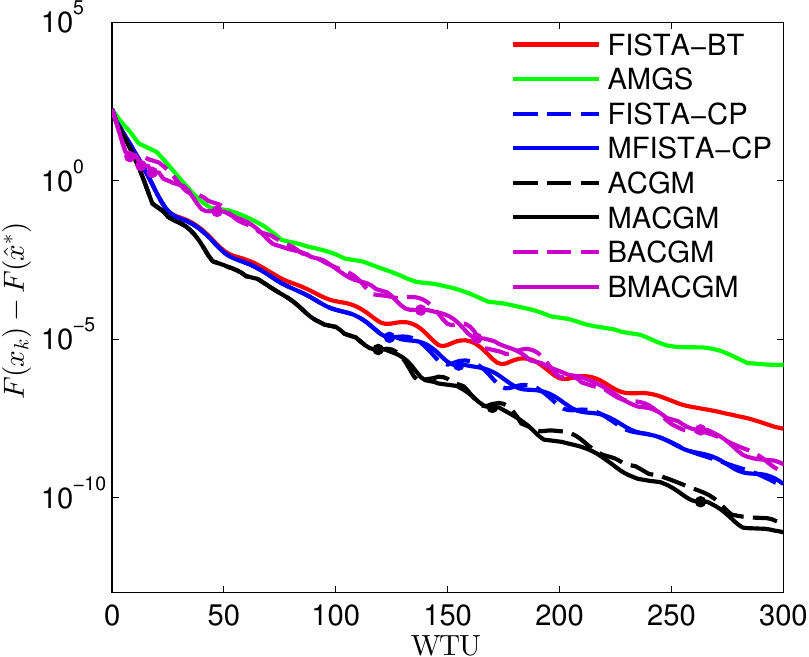}
(d) Computational convergence rates on EN
\end{minipage}

\vspace{6mm}

\begin{minipage}[t]{0.48\linewidth} {\centering
\includegraphics[width=\textwidth]{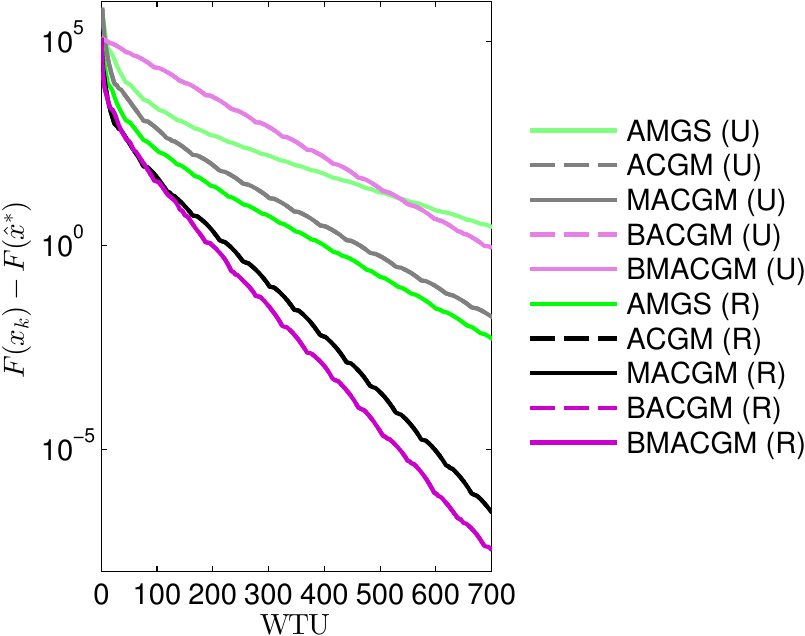}}
(e) Computational convergence rates (R) and upper bounds (U) on RR
\end{minipage}
\hspace{2mm}
\begin{minipage}[t]{0.48\linewidth} {\centering
\includegraphics[width=\textwidth]{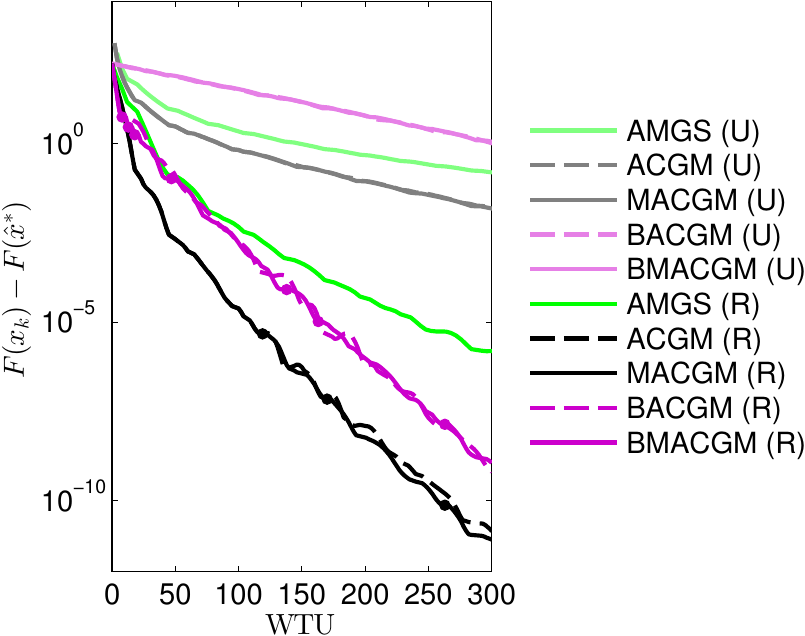}}
(f) Computational convergence rates (R) and upper bounds (U) on EN
\end{minipage}

\vspace{5mm}

\caption{Convergence results of FISTA with backtracking (FISTA-BT), AMGS, FISTA-CP, monotone FISTA-CP (MFISTA-CP), non-monotone ACGM, monotone ACGM (MACGM), border-case non-monotone ACGM (BACGM), and border-case monotone ACGM (BMACGM) on the RR and EN strongly-convex problems. Dots mark iterations preceding overshoots.}
\label{labelaaf}
\end{figure}
}

\newcommand{\cmdaan}{
\begin{figure}[t!] \centering
\begin{minipage}[t]{0.45\linewidth} \centering
\includegraphics[width=\textwidth]{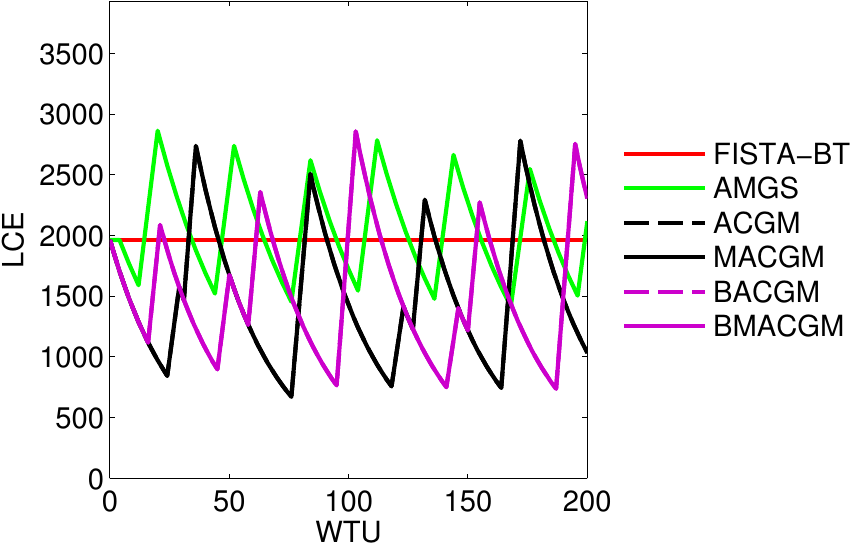}
(a) RR
\end{minipage}
\hspace{10mm}
\begin{minipage}[t]{0.45\linewidth} \centering
\includegraphics[width=\textwidth]{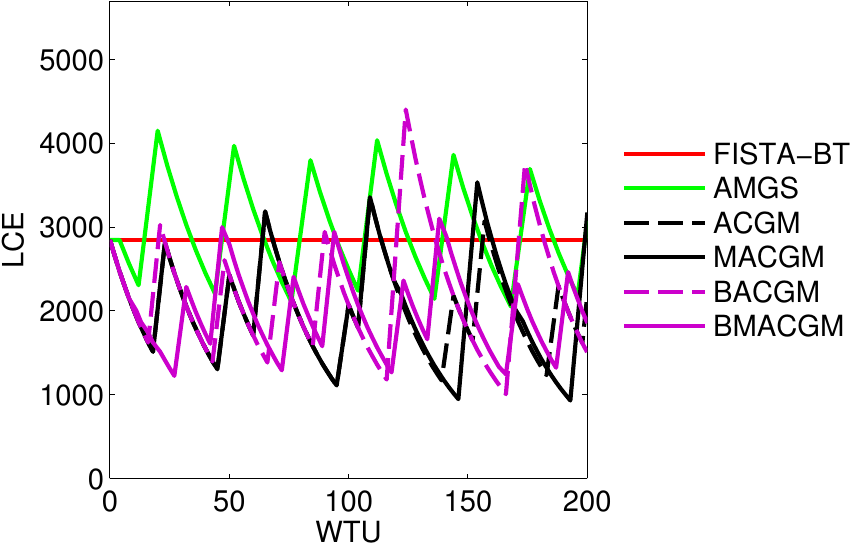}
(b) EN
\end{minipage}
\caption{Line-search method LCE variation on RR and EN}
\label{labelaag}

\vspace{3mm}

\captionof{table}{Average LCEs of line-search methods on RR and EN}
\label{labelaah}
\footnotesize
\begin{tabular}{@{}c@{}c@{}ccccccc@{}} \toprule
Problem & $L_f$ & Iterations & FISTA-BT & AMGS & ACGM & MACGM & BACGM & BMACGM \\ \midrule
RR & 1963.66 & 350 & 1963.66 & 2022.73 & 1473.88 & 1473.88 & 1471.16 & 1471.16 \\
EE & 2846.02 & 150 & 2846.02 & 3023.47 & 2056.68 & 2003.09 & 2093.56 & 1998.12 \\ \bottomrule
\end{tabular}
\end{figure}
}

\newenvironment{pattern}[1][H]
{\floatname{algorithm}{Pattern}
\begin{algorithm}[#1]
}{\end{algorithm}}

\begin{document}

\title{A Generalized Accelerated Composite Gradient Method: Uniting Nesterov's Fast Gradient Method and FISTA}

\author{Mihai I. Florea\thanks
{Department of Signal Processing and Acoustics,
Aalto University, Espoo, Finland
\mbox{(\email{mihai.florea@aalto.fi, sergiy.vorobyov@aalto.fi})}.}
\and
Sergiy A. Vorobyov\samethanks}

\maketitle

\newcommand{\cmdaao}{Generalized Accelerated Composite Gradient Method}
\newcommand{\cmdaap}{Mihai I. Florea and Sergiy A. Vorobyov}
\headers{\cmdaao}{\cmdaap}

\begin{abstract}
Numerous problems in signal processing, statistical inference, computer vision, and machine learning, can be cast as large-scale convex optimization problems. Due to their size, many of these problems can only be addressed by first-order accelerated black-box methods. The most popular among these are the Fast Gradient Method (FGM) and the Fast Iterative Shrinkage Thresholding Algorithm (FISTA). FGM requires that the objective be finite and differentiable with known gradient Lipschitz constant. FISTA is applicable to the more broad class of composite objectives and is equipped with a line-search procedure for estimating the Lipschitz constant. Nonetheless, FISTA cannot increase the step size and is unable to take advantage of strong convexity. FGM and FISTA are very similar in form. Despite this, they appear to have vastly differing convergence analyses. In this work we generalize the previously introduced augmented estimate sequence framework as well as the related notion of the gap sequence. We showcase the flexibility of our tools by constructing a Generalized Accelerated Composite Gradient Method, that unites FGM and FISTA, along with their most popular variants. The Lyapunov property of the generalized gap sequence used in deriving our method implies that both FGM and FISTA are amenable to a Lyapunov analysis, common among optimization algorithms. We further showcase the flexibility of our tools by endowing our method with monotonically decreasing objective function values alongside a versatile line-search procedure. By simultaneously incorporating the strengths of FGM and FISTA, our method is able to surpass both in terms of robustness and usability. We support our findings with simulation results on an extensive benchmark of composite problems. Our experiments show that monotonicity has a stabilizing effect on convergence and challenge the notion present in the literature that for strongly convex objectives, accelerated proximal schemes can be reduced to fixed momentum methods.
\end{abstract}



\section{Introduction}

Numerous large-scale convex optimization problems have recently emerged in a variety of fields, including signal and image processing, statistical inference, computer vision, and machine learning. Often, little is known about the actual structure of the objective function. Therefore, optimization algorithms used in solving such problems can only rely (e.g., by means of callback functions) on specific black-box methods, called oracle functions~\cite{refaab}. The term ``large-scale'' refers to the tractability of certain computational primitives (see also~\cite{refaac}). In the black-box setting, it means that the oracle functions of large-scale problems usually only include scalar functions and operations that resemble first-order derivatives.

Many large-scale applications were rendered practical to address with the advent of Nesterov's Fast Gradient Method (FGM)~\cite{refaad}. The canonical form of FGM, formulated in \cite{refaae}, requires that the objective be differentiable with Lipschitz gradient. Many optimization problems, particularly inverse problems in fields such as sparse signal processing, linear algebra, matrix and tensor completion, and digital imaging (see \cite{refaaf,refaag,refaah,refaai,refaaj} and references therein), have a composite structure. In these composite problems, the objective $F$ is the sum of a function $f$ with Lipschitz gradient (Lipschitz constant $L_f$) and a simple but possibly non-differentiable regularizer $\Psi$. The regularizer $\Psi$ embeds constraints by being infinite outside the feasible set. Often, $L_f$ is not known in advance. The composite problem oracle functions are the scalar $f(\bm{x})$ and $\Psi(\bm{x})$, as well as the gradient $\nabla f(\bm{x})$ and the proximal operator $prox_{\tau \Psi}(\bm{x})$.

To address composite problems, Nesterov has devised an Accelerated Multistep Gradient Scheme (AMGS)~\cite{refaak}. This method updates a Lipschitz constant estimate (LCE) at every iteration using a subprocess commonly referred to in the literature as ``line-search''~\cite{refaai,refaal}. However, the generation of a new iterate (advancement phase of an iteration) and line-search are interdependent and cannot be executed in parallel. Moreover, AMGS utilizes only the gradient-type oracle functions $\nabla f(\bm{x})$ and $prox_{\tau \Psi}(\bm{x})$. In many applications, including compressed sensing (e.g., LASSO~\cite{refaam}) and many classification tasks (e.g., $l_1$-regularized logistic regression), the evaluation of $\nabla f(\bm{x})$ is more computationally expensive than $f(\bm{x})$. An alternative to AMGS that uses $f(\bm{x})$ calls in line-search has been proposed by Beck and Teboulle in the form of the Fast Iterative Shrinkage-Thresholding Algorithm (FISTA)~\cite{refaal}.\footnote{FISTA has very similar functionality and convergence guarantees to those of the first variant of FGM introduced in \cite{refaad}. For clarity, we denote the latter method also as FISTA, whereas we use FGM to designate the canonical form developed in \cite{refaae}.} FISTA also benefits from having line-search decoupled from advancement. However, FISTA is unable to decrease the LCE at run-time and does not take into account strong convexity.

A strongly convex extension of FISTA, which we designate as FISTA Chambolle-Pock (FISTA-CP), has been recently introduced in~\cite{refaah}, but without line-search. Two works~\cite{refaan, refaao} propose a similar method that accommodates only feasible start, assumes $\mu_{\Psi} = 0$, and does not feature line-search nor monotonically decreasing objective function values (monotonicity). In the context of the composite problem class, this method is merely a restricted version of FISTA-CP. Consequently, we also treat it as FISTA-CP. Another study~\cite{refaap} proposes a variant of FISTA-CP equipped with a line-search procedure capable of decreasing the LCE (a fully adaptive line-search strategy). However, this strongly convex Accelerated Proximal Gradient (scAPG) method is guaranteed to converge only if $f$ is strongly convex and if the starting point is feasible.

Apart from the above methods, which are applicable to a wide range of problems and feature convergence guarantees, the literature abounds with algorithms that either employ heuristics or are highly specialized. For instance, restart techniques have been developed to address strong convexity. However, the adaptive restart proposed in \cite{refaaq} lacks rigorous convergence guarantees for non-quadratic objectives and the periodic restart introduced in \cite{refaad} requires prior knowledge of the Lipschitz constant. A very interesting study~\cite{refaar} proposes an over-relaxed variant of FISTA, a first-order method that appears to define a distinct category of algorithms. However, the study fails to argue whether the convergence guarantees are superior to those of FISTA. Despite being tested only on instances of the LASSO problem~\cite{refaam}, the performance benefit is unclear even in this restricted scenario. Numerous other methods successfully exploit additional problem structure, such as sparsity of optimal points (e.g.~\cite{refaas,refaat,refaau,refaav}), however they cannot be applied to the entire composite problem class.

FISTA (when the objective is non-strongly convex), FISTA-CP, and scAPG (when $f$ is strongly-convex) are identical in form to certain variants of FGM. Therefore, FISTA, FISTA-CP, and scAPG can be viewed as extensions of FGM for composite objectives. However, whereas FGM was derived using the estimate sequence~\cite{refaae}, FISTA, FISTA-CP (including the variant found in~\cite{refaan, refaao}), and scAPG were proposed without derivation. No details have been provided on how the update rules of each method were constructed, apart from their obvious validation in the corresponding convergence analysis. Consequently, new features of FGM cannot be directly incorporated into FISTA, FISTA-CP, and scAPG. For instance, Nesterov has proposed in \cite{refaaw} a line-search variant of FGM with an implicit derivation based on the estimate sequence, albeit only for non-strongly convex objectives. Had FISTA been derived in this way, a fully adaptive line-search variant of FISTA could simply be obtained from \cite{refaaw}. Instead, a sophisticated fully adaptive line-search extension was proposed in~\cite{refaax}, with a technical derivation based on the mathematical constructs of \cite{refaal}. On the other hand, through partial adoption of the estimate sequence, i.e., relating FISTA to the ``constant step scheme I'' variant of FGM in~\cite{refaae} and AMGS, we have arrived at a similar but simpler fully adaptive line-search scheme for FISTA~\cite{refaay}, but again only in the non-strongly convex scenario.

In~\cite{refaaz}, we have introduced the \emph{augmented estimate sequence} framework and used it to derive the Accelerated Composite Gradient Method (ACGM), which incorporates by design a fully adaptive line-search procedure. ACGM has the convergence guarantees of FGM, the best among primal first-order methods, while being as broadly applicable as AMGS. In addition, FISTA-CP and FISTA, along with the fully adaptive line-search extensions in \cite{refaax} and \cite{refaay}, are particular cases of ACGM~\cite{refaaz}. However, to accommodate infeasible start, we have imposed restrictions on the input parameters. Algorithms such as the ``constant step scheme III'' variant of FGM~\cite{refaae} and scAPG~\cite{refaap} are guaranteed to converge only when the starting point is feasible, and thus do not correspond to any instance of ACGM in~\cite{refaaz}.

\subsection{Contributions}

\begin{itemize}

\item In this paper, by lifting the restrictions on the input parameters imposed by infeasible start, we generalize the augmented estimate sequence framework and derive a generalization of ACGM that encompasses FGM and FISTA, along with their variants.

\item We further showcase the flexibility and power of the augmented estimate sequence framework by endowing ACGM with monotonicity alongside its adaptive line-search procedure. Monotonicity is a desirable property, particularly when dealing with proximal operators that lack a closed form expression or other kinds of inexact oracles~\cite{refaah,refaba}. Even when dealing with exact oracles, monotonicity leads to a more stable and predictable convergence rate.

\item Our generalized ACGM, with its additional features, including the fully adaptive estimation of the Lipschitz constant, is therefore superior to FGM and FISTA in terms of flexibility and usability.

\item We further generalize the definition of the wall-clock time unit (previously introduced in~\cite{refaaz}), to account for arbitrary oracle function costs. We analyze the per-iteration complexity of our method in this generalized context, taking advantage of parallelization when possible.

\item We support our theoretical findings with simulation results on an extended benchmark of signal processing and machine learning problems.

\end{itemize}

\subsection{Assumptions and notation} \label{labelaai}

We consider the composite optimization problems of the form
\begin{equation}
\displaystyle \min_{\bm{x}\in \mathbb{R}^n} F(\bm{x}) \operatorname{\ \overset{def}{=} \ } f(\bm{x}) + \Psi(\bm{x}),
\end{equation}
where $\bm{x}$ is a vector of $n$ optimization variables (bold type indicates a vector in $\mathbb{R}^n$), and $F$ is the objective function. The constituents of the objective $F$ are the convex differentiable function \mbox{$f : \mathbb{R}^n \rightarrow \mathbb{R}$} and the convex lower semicontinuous regularizer function \mbox{$\Psi : \mathbb{R}^n \rightarrow \mathbb{R} \cup \{\infty\}$}. Function $f$ has Lipschitz gradient (Lipschitz constant $L_f > 0$) and a strong convexity parameter $\mu_f \geq 0$. Note that the existence of a strong convexity parameter does not imply the strong convexity property since any convex function has a strong convexity parameter of zero. Regularizer $\Psi$ has a strong convexity parameter $\mu_{\Psi} \geq 0$, entailing that objective $F$ has a strong convexity parameter $\mu = \mu_f + \mu_{\Psi}$. Constraints are enforced by making $\Psi$ infinite outside the feasible set, which is closed and convex.

Apart from the above properties, nothing is assumed known about functions $f$ and $\Psi$, which can only be accessed in a \emph{black-box} fashion \cite{refaab} by querying oracle functions $f(\bm{x})$, $\nabla f(\bm{x})$, $\Psi(\bm{x} )$, and $\mbox{prox}_{\tau \Psi}(\bm{x})$, with arguments $\bm{x} \in \mathbb{R}^n$ and $\tau > 0$. The proximal operator $\mbox{prox}_{\tau \Psi}(\bm{x})$ is given by
\begin{equation}
\mbox{prox}_{\tau \Psi}(\bm{x}) \operatorname{\ \overset{def}{=} \ } \displaystyle \argmin_{\bm{z} \in \mathbb{R}^n}\left(\Psi(\bm{z}) + \frac{1}{2\tau}\|\bm{z} - \bm{x}\|_{2}^{2}\right),
\end{equation}
for all $\bm{x} \in \mathbb{R}^n$ and $\tau > 0$, with $\| . \|_2$ designating the standard Euclidean norm in $\mathbb{R}^n$.

Central to our derivation are \emph{generalized parabolae}, quadratic functions whose Hessians are multiples of the identify matrix. We refer to the strongly convex ones simply as \emph{parabolae}. Any parabola $\psi : \mathbb{R}^n \rightarrow \mathbb{R}$ can be written in canonical form as
\begin{equation}
\psi(\bm{x}) \operatorname{\ \overset{def}{=} \ } \psi^* + \frac{\gamma}{2} \|\bm{x} - \bm{v}\|_2^2, \quad\bm{x} \in \mathbb{R}^n,
\end{equation}
where $\gamma > 0$ denotes the curvature, $\bm{v}$ is the vertex, and $\psi^*$ is the optimum value.

For conciseness, we introduce the generalized parabola expression $Q_{\lambda, \bm{y}}(\bm{x})$ for all $\bm{x}, \bm{y} \in \mathbb{R}^n$ and $\lambda \geq 0$ as
\begin{equation} \label{labelaaj}
Q_{\lambda, \bm{y}}(\bm{x}) \operatorname{\ \overset{def}{=} \ } f(\bm{y}) + \langle \nabla f(\bm{y}), \bm{x} - \bm{y} \rangle + \frac{\lambda}{2} \|\bm{x} - \bm{y} \|_2^2.
\end{equation}
The proximal gradient operator $T_{L}(\bm{y})$ \cite{refaak} can be expressed succinctly using \eqref{labelaaj} as
\begin{equation} \label{labelaak}
T_{L}(\bm{y}) \operatorname{\ \overset{def}{=} \ } \displaystyle \argmin_{\bm{x} \in \mathbb{R}^n} \left( Q_{L, \bm{y}}(\bm{x}) + \Psi(\bm{x}) \right), \quad \bm{y} \in \mathbb{R}^n,
\end{equation}
where $L > 0$ is a parameter corresponding to the inverse of the step size. Operator $T_{L}(\bm{y})$ can be evaluated in terms of oracle functions as
\begin{equation}
T_{L}(\bm{y}) = \mbox{prox}_{\frac{1}{L} \Psi} \left(\bm{y} - \frac{1}{L} \nabla f(\bm{y}) \right), \quad \bm{y} \in \mathbb{R}^n, \quad L > 0.
\end{equation}

\section{Generalizing ACGM}

In our previous formulation of the augmented estimate sequence framework and subsequent derivation of the Accelerated Composite Gradient Method (ACGM) in~\cite{refaaz}, we have chosen to accommodate infeasible start by default. In the sequel, we argue why this assumption was necessary and provide a generalized framework, wherein this assumption is only a particular case.

\subsection{FGM estimate sequence} \label{labelaal}

We begin the generalization of the estimate sequence by noting that a rate of convergence can only be obtained on function values and not iterates. It may be argued that, when the objective function is strongly convex, many first-order schemes, including the non-accelerated fixed-point methods, guarantee the linear convergence of the iterates to the optimal point. However, when the problem is non-strongly convex, the optimization landscape may contain a high-dimensional subspace of very low curvature in the vicinity of the set of optimal points. In this case, a convergence rate of iterates remains a difficult open problem~\cite{refabb}. For instance, Nesterov has provided in \cite{refaae} an ill-conditioned quadratic problem where a convergence rate for iterates cannot be formulated for all first-order schemes of a certain structure.

Hence, we choose to measure convergence using the image space distance (ISD), which is the distance between the objective values at iterates and the optimal value. The decrease rate of an upper bound on the ISD gives the convergence guarantee (provable convergence rate). The generalization of the estimate sequence framework follows naturally from the formulation of such guarantees. Specifically, we interpret the image space distance upper bound (ISDUB), provided by Nesterov for FGM in~\cite{refaae}, for all $k \geq 0$ as
\begin{equation} \label{labelaam}
A_k (F(\bm{x}_k) - F(\bm{x}^*)) \leq A_0( F(\bm{x}_0) - F(\bm{x}^*)) + \frac{\gamma_0}{2} \|\bm{x}_0 - \bm{x}^* \|_2^2.
\end{equation}
Point $\bm{x}^*$ can be any optimal point. However, we consider it fixed throughout this work. The convergence guarantee is given by the sequence $\{A_k\}_{k \geq 0}$ with $A_0 \geq 0$. Each iteration may not be able to improve the objective function value, but at least it must improve the convergence guarantee, meaning that $\{A_k\}_{k \geq 0}$ must be an increasing sequence. From $A_0 \geq 0$ it follows that $A_k > 0$ for all $k \geq 1$.

The right-hand side of~\eqref{labelaam} is a weighted sum between the initial ISD and the corresponding domain space term (DST), with weights given by $A_0$ and $\gamma_0$, respectively. In the derivation of FGM, the weights are constrained as $A_0 > 0$ and $\gamma_0 \geq A_0 \mu$ whereas in the original ACGM \cite{refaaz} we have enforced $A_0 = 0$ and $\gamma_0 \leq 1$.
Given that our current aim is to provide a generic framework, we impose no restrictions on the weights, apart from $A_0 \geq 0$ and $\gamma_0 > 0$. The former restriction follows from the convexity of $F$ while the latter is required by the estimate sequence, along with its augmented variant, as we shall demonstrate in the sequel.

The ISDUB expression can be rearranged to take the form
\begin{equation} \label{labelaan}
A_k F(\bm{x}_k) \leq H_k,
\end{equation}
where
\begin{equation}
H_k \operatorname{\ \overset{def}{=} \ } (A_k - A_0) F(\bm{x}^*) + A_0 F(\bm{x}_0) + \frac{\gamma_0}{2} \|\bm{x}_0 - \bm{x}^* \|_2^2
\end{equation}
is the highest upper bound that can be placed on weighted objective values $A_k F(\bm{x}_k)$ to satisfy \eqref{labelaam}.

The value of $H_k$ depends on the optimal value $F(\bm{x}^*)$, which is an unknown quantity. The estimate sequence provides a computable, albeit more stringent, replacement for $H_k$. It is obtained as follows. The convexity of the objective implies the existence of a sequence $\{W_k\}_{k \geq 1}$ of convex global lower bounds on $F$, written as
\begin{equation} \label{labelaao}
F(\bm{x}) \geq W_k(\bm{x}) \mbox{, } \quad \bm{x} \in \mathbb{R}^n , \quad k \geq 1.
\end{equation}

For instance, when querying $\nabla f(\bm{x})$ at a test point $\bm{y}_{k}$, we obtain a global lower bound on $f$ as $Q_{\mu_f, \bm{y}_k}(\bm{x})$ for all $\bm{x} \in \mathbb{R}^n$. Consequently $Q_{\mu_f, \bm{y}_k} + \Psi$ is a global lower bound of $F$. The lower bounds obtained from test points up to and including $\bm{y}_k$ can be combined to produce $W_k$. The drawback of incorporating $\Psi$ in the lower bounds lies in the introduction of additional $\mbox{prox}_{\tau \Psi}(\bm{x})$ calls, as in~\cite{refaak}. In Subsection~\ref{labelaay} we will formulate more computationally efficient bounds. Notwithstanding, the framework we introduce here does not assume any structure on $W_k$ apart from convexity.

By substituting the optimal value terms $F(\bm{x}^*)$ in \eqref{labelaan} with $W_k(\bm{x}^*)$, we obtain $\mathcal{H}_k$, a lower bound on $H_k$, given by
\begin{equation}
\mathcal{H}_k \operatorname{\ \overset{def}{=} \ } (A_k - A_0) W_k(\bm{x}^*) + A_0 F(\bm{x}_0) + \frac{\gamma_0}{2} \| \bm{x}^* - \bm{x}_0\|_2^2,
\end{equation}
for all $k \geq 0$. This still depends on $\bm{x}^*$. However, $\mathcal{H}_k$ can be viewed as the value of an \emph{estimate function}, taken at an optimal point $\bm{x}^*$. The estimate functions $\psi_k$ are defined as functional extensions of $\mathcal{H}_k$, namely
\begin{equation} \label{labelaap}
\psi_k(\bm{x}) \operatorname{\ \overset{def}{=} \ } (A_k - A_0) W_k(\bm{x}) + A_0 F(\bm{x}_0) + \frac{\gamma_0}{2} \|\bm{x} - \bm{x}_0 \|_2^2,
\end{equation}
for all $\bm{x} \in \mathbb{R}^n$ and $k \geq 0$.
Note that the first estimate function $\psi_0$ does not contain a lower bound term. Therefore, it is not necessary to define $W_0$. The collection of estimate functions $\{\psi_k\}_{k \geq 0}$, is referred to as the \emph{estimate sequence}.

The estimate function optimum value, given by
\begin{equation} \label{labelaaq}
\psi_k^* \operatorname{\ \overset{def}{=} \ } \displaystyle \min_{\bm{x} \in \mathbb{R}^n} \psi_k(\bm{x}) , \quad k \geq 0,
\end{equation}
is guaranteed to be lower than $\mathcal{H}_k$, since
\begin{equation} \label{labelaar}
\psi_k^* = \displaystyle \min_{\bm{x} \in \mathbb{R}^n} \psi_k(\bm{x}) \leq \psi_k(\bm{x}^*) = \mathcal{H}_k , \quad k \geq 0.
\end{equation}
As such, $\psi_k^*$ provides the sought after computable replacement of $H_k$. Note that if the lower bounds $W_k$ are linear, the estimate functions are generalized parabolae with the curvature given by $\gamma_0$. In this case, the existence of $\psi_k^*$ is conditioned by $\gamma_0 > 0$, explaining the assumption made in~\eqref{labelaam}.

Thus, it suffices to maintain the estimate sequence property, given by
\begin{equation} \label{labelaas}
A_k F(\bm{x}_k) \leq \psi_k^* , \quad k \geq 0,
\end{equation}
to satisfy the ISDUB expression~\eqref{labelaam}. The proof follows from the above definitions as
\begin{equation}
A_k F(\bm{x}_k) \overset{\eqref{labelaas}}{\leq}
\psi_k^* \overset{\eqref{labelaar}}{\leq}
\psi_k(\bm{x}^*) =
\mathcal{H}_k \overset{\eqref{labelaao}}{\leq} H_k , \quad k \geq 0.
\end{equation}

\subsection{Generalizing the augmented estimate sequence} \label{labelaat}

The interval between the maintained upper bound $\psi_k^*$ and the highest allowable bound $H_k$ contains $\mathcal{H}_k$. This allows us to produce a relaxation of the estimate sequence by forcibly \emph{closing the gap} between $\mathcal{H}_k$ and $H_k$. Namely, we define the augmented estimate functions as
\begin{equation} \label{labelaau}
\psi'_k (\bm{x}) \operatorname{\ \overset{def}{=} \ } \psi_k(\bm{x}) + H_k - \mathcal{H}_k , \quad k \geq 0,
\end{equation}
with $\{\psi'_k(\bm{x})\}_{k \geq 0}$ being the augmented estimate sequence. We expand definition~\eqref{labelaau} as
\begin{equation} \label{labelaav}
\psi'_k (\bm{x}) = \psi_k(\bm{x}) + (A_k - A_0) (F(\bm{x}^*) - W_k(\bm{x}^*)) , \quad k \geq 0.
\end{equation}
The augmented estimate sequence property is formulated as
\begin{equation} \label{labelaaw}
A_k F(\bm{x}_k) \leq \psi'^*_k , \quad k \geq 0,
\end{equation}
where the augmented estimate function optimum value $\psi'^*_k$ is given by
\begin{equation} \label{labelaax}
\psi_k'^* \operatorname{\ \overset{def}{=} \ } \displaystyle \min_{\bm{x} \in \mathbb{R}^n} \psi'_k(\bm{x}) , \quad k \geq 0.
\end{equation}

\subsection{Quadratic lower bounds} \label{labelaay}

As we have seen in Subsection~\ref{labelaal}, the structure of the estimate sequence and, subsequently, the augmented estimate sequence, depends on the global lower bounds $W_k$. These can be obtained by combining simple lower bounds, each produced at each iteration $k$ using test point $\bm{y}_k$ and denoted as $w_k$. In the following, we propose simple lower bounds $w_k$ that can be obtained in a computationally efficient manner, each requiring only one call to $T_L(\bm{x})$.

\begin{lemma} \label{labelaaz}
If at every iteration $k \geq 0$, the auxiliary points and LCEs obey the descent condition, stated as
\begin{equation} \label{labelaba}
f(\bm{z}_{k + 1}) \leq Q_{L_{k + 1}, \bm{y}_{k + 1}} (\bm{z}_{k + 1}) , \quad k \geq 0,
\end{equation}
where
\begin{equation} \label{labelabb}
\bm{z}_{k + 1} \operatorname{\ \overset{def}{=} \ } T_{L_{k + 1}}(\bm{y}_{k + 1}),
\end{equation}
then the objective is lower bounded for all $\bm{x} \in \mathbb{R}^n$ as
\begin{equation}
\begin{gathered}
F(\bm{x}) \geq w_{k + 1}(\bm{x}) \operatorname{\ \overset{def}{=} \ } F(\bm{z}_{k + 1}) + \frac{L_{k + 1} + \mu_{\Psi}}{2} \| \bm{z}_{k + 1} - \bm{y}_{k + 1} \|_2^2 \\
+ (L_{k + 1} + \mu_{\Psi}) \langle \bm{y}_{k + 1} - \bm{z}_{k + 1}, \bm{x} - \bm{y}_{k + 1} \rangle + \frac{\mu}{2} \| \bm{x} - \bm{y}_{k + 1} \|_2^2.
\end{gathered}
\end{equation}
\end{lemma}
\begin{proof}
By expanding \eqref{labelabb} we obtain
\begin{equation}
\bm{z}_{k + 1}= \displaystyle \argmin_{\bm{x} \in \mathbb{R}^n} \left( f(\bm{y}_{k + 1}) + \langle \nabla f(\bm{y}_{k + 1}), \bm{x} - \bm{y}_{k + 1} \rangle + \frac{L_{k + 1}}{2} \| \bm{x} - \bm{y}_{k + 1} \|_2^2 + \Psi(\bm{x}) \right).
\end{equation}
This subproblem is unconstrained and has a strongly convex and subdifferentiable objective. The first-order optimality condition is therefore equivalent to
\begin{equation} \label{labelabc}
\bm{\xi} = L_{k + 1} (\bm{y}_{k + 1} - \bm{z}_{k + 1}) - \nabla f(\bm{y}_{k + 1}),
\end{equation}
where $\bm{\xi}$ is a specific subgradient of function $\Psi$ at point $\bm{z}_{k + 1}$.

Moreover, from the definition of the strong convexity parameters of $f$ and $\Psi$ we have that
\begin{align}
f(\bm{x}) &\geq f(\bm{y}_{k + 1}) + \langle \nabla f(\bm{y}_{k + 1}), \bm{x} - \bm{y}_{k + 1} \rangle + \frac{\mu_f}{2} \| \bm{x} - \bm{y}_{k + 1} \|_2^2 , \label{labelabd}\\
\Psi(\bm{x}) &\geq \Psi(\bm{z}_{k + 1}) + \langle \bm{\xi}, \bm{x} - \bm{z}_{k + 1} \rangle + \frac{\mu_{\Psi}}{2} \| \bm{x} - \bm{z}_{k + 1} \|_2^2 . \label{labelabe}
\end{align}
By adding together \eqref{labelaba}, \eqref{labelabd}, and \eqref{labelabe}, canceling and rearranging terms we obtain
\begin{equation} \label{labelabf}
\begin{aligned}
F(\bm{x}) \geq \ &F(\bm{z}_{k + 1}) + \langle \nabla f(\bm{y}_{k + 1}) + \bm{\xi}, \bm{x} - \bm{z}_{k + 1} \rangle + \frac{\mu_f}{2} \| \bm{x} - \bm{y}_{k + 1} \|_2^2 \\
& + \frac{\mu_{\Psi}}{2} \| \bm{x} - \bm{z}_{k + 1} \|_2^2 - \frac{L_{k + 1}}{2} \| \bm{y}_{k + 1} - \bm{z}_{k + 1} \|_2^2.
\end{aligned}
\end{equation}
Further applying \eqref{labelabc} to \eqref{labelabf} and rearranging terms gives the desired result.
\end{proof}

Once we have obtained the simple lower bounds, we need to determine a means of combining them. We propose the convex combination given by
\begin{equation} \label{labelabg}
W_{k + 1}(\bm{x}) = \left( \displaystyle \sum_{i = 1}^{k + 1}{a_i} \right)^{-1} \left(\displaystyle \sum_{i = 1}^{k + 1}{a_i w_i(\bm{x})} \right) , \quad k \geq 0,
\end{equation}
where
\begin{equation} \label{labelabh}
A_{k + 1} = A_k + a_{k + 1}, \quad k \geq 0.
\end{equation}
We have assumed that $\{A_k\}_{k \geq 0}$ is increasing which implies that the weights are strictly positive, namely
\begin{equation} \label{labelabi}
a_{k + 1} > 0 , \quad k \geq 0.
\end{equation}

From definition~\eqref{labelaap}, we have that the initial estimate function is a parabola, given by
\begin{equation} \label{labelabj}
\psi_0(\bm{x}) = A_0 F(\bm{x}_0) + \frac{\gamma_0}{2} \|\bm{x} - \bm{x}_0 \|_2^2 ,\quad \bm{x} \in \mathbb{R}^n.
\end{equation}
Applying \eqref{labelabg} and \eqref{labelabh} to \eqref{labelaap} we get simple recursion rule in the form of
\begin{equation} \label{labelabk}
\psi_{k + 1}(\bm{x}) = \psi_k(\bm{x}) + a_{k + 1} w_{k + 1}(\bm{x}),\quad \bm{x} \in \mathbb{R}^n, \quad k \geq 0.
\end{equation}

\begin{lemma} \label{labelabl}
When the simple lower bounds are given by Lemma~\ref{labelaaz} and are combined as in \eqref{labelabg}, then the estimate functions and the augmented estimate functions at every iteration are parabolic, and thus can be written in canonical form for all $\bm{x} \in \mathbb{R}^n$ and $k \geq 0$ as
\begin{align}
\psi _k(\bm{x}) &= \psi_k ^* + \frac{\gamma_k}{2} \|\bm{x} - \bm{v}_k\|_2^2 , \label{labelabm} \\
\psi'_k(\bm{x}) &= \psi_k'^* + \frac{\gamma_k}{2} \|\bm{x} - \bm{v}_k\|_2^2 , \label{labelabn}
\end{align}
with $\psi_k^*$ and $\psi_k'^*$ given by \eqref{labelaaq} and \eqref{labelaax}, respectively, and
\begin{equation} \label{labelabo}
\bm{v}_k \operatorname{\ \overset{def}{=} \ } \argmin_{\bm{x} \in \mathbb{R}^n} \psi _k(\bm{x}).
\end{equation}
Moreover, the curvatures and vertices obey the following recursion rules
\begin{gather}
\gamma_{k + 1} = \gamma_k + a_{k + 1} \mu, \label{labelabp} \\
\bm{v}_{k + 1} = \frac{1}{\gamma_{k + 1}} (\gamma_k \bm{v}_k + a_{k + 1} (L_{k + 1} + \mu_{\Psi}) \bm{z}_{k + 1} - a_{k + 1} (L_{k + 1} - \mu_f) \bm{y}_{k + 1} ). \label{labelabq}
\end{gather}
\end{lemma}
\begin{proof}
The Hessian of $\psi_0(\bm{x})$ is given by $\gamma_0 I_n$ (where $I_n$ is the identity matrix of size $n$) with $\gamma_0 > 0$ while the Hessian of $w_{k + 1}(\bm{x})$ is $\mu I_n$. From the recursion rule in \eqref{labelabk}, it follows by simple induction that all estimate functions are parabolae, and as such, can be written as in \eqref{labelabm} with \eqref{labelaaq} and \eqref{labelabo}. The relation in \eqref{labelaav} means that the augmented estimate functions differ from the estimate function by a constant term. It follows that \eqref{labelabn} with \eqref{labelaax} holds.

Expanding the recursion rule in \eqref{labelabk} using \eqref{labelabm} and taking the gradient we obtain
\begin{equation} \label{labelabr}
\begin{gathered}
\gamma_{k + 1} (\bm{x} - \bm{v}_{k + 1}) = \gamma_{k} (\bm{x} - \bm{v}_{k}) + a_{k + 1}(L_{k + 1} + \mu_{\Psi})(\bm{y}_{k + 1} - \bm{z}_{k + 1}) \\
\quad + \ a_{k + 1} \mu (\bm{x} - \bm{y}_{k + 1}), \quad \bm{x} \in \mathbb{R}^n.
\end{gathered}
\end{equation}
Matching the first-order coefficients and constant terms on both sides of \eqref{labelabr} gives \eqref{labelabp} and \eqref{labelabq}, respectively.
\end{proof}
By iterating \eqref{labelabp}, we obtain a closed form expression of the estimate function curvature as
\begin{equation} \label{labelabs}
\gamma_k = \gamma_0 + \left( \displaystyle \sum_{i = 1}^k a_i \right) \mu = \gamma_0 - A_0 \mu + A_k \mu , \quad k \geq 0.
\end{equation}

\subsection{Generalizing the gap sequence}

We introduce the augmented estimate sequence gaps $\{\Gamma_k\}_{k \geq 0}$ and the gap sequence $\{ \Delta_k \}_{k \geq 0}$, respectively, as
\begin{gather}
\Gamma_k \operatorname{\ \overset{def}{=} \ } A_k F(\bm{x}_k) - \psi'^*_k , \quad k \geq 0, \label{labelabt} \\
\Delta_k \operatorname{\ \overset{def}{=} \ } A_k (F(\bm{x}_k) - F(\bm{x}^*)) + \frac{\gamma_k}{2}\| \bm{v}_k - \bm{x}^* \|_2^2 , \quad k \geq 0 . \label{labelabu}
\end{gather}

\begin{proposition} \label{labelabv}
When the estimate functions are parabolic, as in Lemma~\ref{labelabl}, the augmented estimate sequence gaps can be expressed more succinctly as
\begin{equation} \label{labelabw}
\Gamma_k = \Delta_k - \Delta_0 , \quad k \geq 0.
\end{equation}
\end{proposition}
\begin{proof}
From \eqref{labelaav}, using \eqref{labelabm} and \eqref{labelabn} in Lemma~\ref{labelabl}, it follows that
\begin{equation} \label{labelabx}
\psi'^*_k = \psi^*_k + (A_k - A_0) (F(\bm{x}^*) - W_k(\bm{x}^*)).
\end{equation}
Therefore, the gap $\Gamma_k$ can be expressed as
\begin{align}
\Gamma_k
\overset{\eqref{labelabx}}{=}
&A_k (F(\bm{x}_k) - F(\bm{x}^*)) + (A_k - A_0) W_k(\bm{x}^*) + A_0 F(\bm{x}^*) - \psi_k^* \nonumber \\
\overset{\eqref{labelabm}}{=}
&A_k (F(\bm{x}_k) - F(\bm{x}^*)) + \frac{\gamma_k}{2} \|\bm{x}^* - \bm{v_k}\|_2^2 - \nonumber \\
& - \psi_k(\bm{x}^*) + A_0 F(\bm{x}^*) + (A_k - A_0) W_k(\bm{x}^*) \nonumber \\
\overset{\eqref{labelaap}}{=}
&A_k (F(\bm{x}_k) - F(\bm{x}^*)) + \frac{\gamma_k}{2} \|\bm{v_k} - \bm{x}^*\|_2^2 - \nonumber \\
& - A_0 (F(\bm{x}_0) - F(\bm{x}^*)) - \frac{\gamma_0}{2} \|\bm{x_0} - \bm{x}^*\|_2^2 \nonumber \\
\overset{\eqref{labelabu}}{=} &\Delta_k - \Delta_0 , \quad k \geq 0 . \nonumber
\end{align}
\end{proof}
Note that \eqref{labelaav}, \eqref{labelabj}, and Lemma~\ref{labelabl} imply that $\psi_0'^* = \psi_0'(\bm{x}_0) = \psi_0(\bm{x}_0) = A_0 F(\bm{x}_0)$ which, together with \eqref{labelabt}, guarantees that the initial augmented estimate sequence gap $\Gamma_0$ is zero, regardless of the structure of the lower bounds. Therefore, a sufficient condition for the preservation of the augmented estimate sequence property \eqref{labelaaw} as the algorithm progresses is that the augmented estimate sequence gaps do not increase. Proposition~\ref{labelabv} states that when the lower bounds are generalized parabolae, the sufficient condition for \eqref{labelaaw} becomes
\begin{equation} \label{labelaby}
\Delta_{k + 1} \leq \Delta_k , \quad k \geq 0.
\end{equation}
Therefore, \eqref{labelaaw} and the subsequent the convergence of the algorithm are implied by the \emph{Lyapunov} (non-increasing) property of the gap sequence. Such Lyapunov functions have been widely employed in the theoretical analysis of optimization schemes (e.g., \cite{refabc,refabd,refabe}). Moreover, the simple structure and generic nature of the gap sequence further justifies the augmentation of the estimate sequence.

\subsection{A design pattern for first-order accelerated methods} \label{labelabz}

The augmented estimate sequence property places an upper bound on $F(\bm{x}_{k + 1})$, which must be satisfied before the generation of $\bm{x}_{k + 1}$. The black-box nature of the objective constitutes an obstacle in this respect and it is therefore advantageous to employ a surrogate upper bound on the objective. Numerous algorithms employ this technique, and are often referred to as majorization minimization (MM) algorithms (see, e.g., \cite{refabf,refabg} and references therein). New iterates are generated as
\begin{equation} \label{labelaca}
\bm{x}_{k + 1} = \argmin_{\bm{x} \in \mathbb{R}^n} u_{k + 1} (\bm{x}) , \quad k \geq 0,
\end{equation}
where $u_{k + 1}$ is a \emph{local} upper bound on the objective at the new iterate, namely
\begin{equation} \label{labelacb}
F(\bm{x}_{k + 1}) \leq u_{k + 1}(\bm{x}_{k + 1}) , \quad k \geq 0.
\end{equation}
Combining \eqref{labelabh}, \eqref{labelabk}, and \eqref{labelaca} yields the algorithm design pattern outlined in Pattern~\ref{labelacc}. Note that Nesterov's FGM and AMGS adhere to Pattern~\ref{labelacc}. This pattern will form the scaffolding of our generalized ACGM.

Pattern~\ref{labelacc} takes as input the starting point $\bm{x}_0$, the curvature $\gamma_0$, the initial guarantee $A_0$ and, if the Lipschitz constant is not known in advance, an initial LCE $L_0 > 0$. In \cmdaai{labelacc}{labelaci}, the main iterate is given by the minimum of the \emph{local} upper bound $u_{k + 1}$. Alongside the main iterate, the algorithm maintains an estimate function $\psi_{k + 1}$, obtained from the previous one by adding a simple \emph{global} lower bound $w_{k + 1}$ weighted by \mbox{$a_{k + 1} > 0$} (\cmdaai{labelacc}{labelacj}). Function $w_{k + 1}$ need not be given by Lemma~\ref{labelaaz}, but we assume that it is uniquely determined by an auxiliary point $\bm{y}_{k + 1}$. The current LCE $L_{k + 1}$, weight $a_{k + 1}$, and auxiliary point $\bm{y}_{k + 1}$ are computed using algorithm specific methods $\mathcal{S}$, $\mathcal{F}_a$, and $\mathcal{F}_y$, respectively (lines \ref{labelace}, \ref{labelacf}, and \ref{labelacg} of Pattern~\ref{labelacc}). These methods take as parameters the state of the algorithm, given by current values of the main iterate, LCE, weight, and estimate function.

\begin{pattern}[H]
\caption{A design pattern for first-order accelerated methods}
\label{labelacc}
\begin{algorithmic}[1]
\STATE $\psi_0(\bm{x}) = A_0 F(\bm{x}_0) + \frac{\gamma_0}{2} \|\bm{x} - \bm{x}_0 \|_2^2 $ \label{labelacd}
\FOR{$k = 0,\ldots{},K-1$}
\STATE $ L_{k + 1} = \mathcal{S}(\bm{x}_k, \psi_k, A_k, L_k)$ \hfill ``line-search'' \label{labelace}
\STATE $ a_{k + 1} = \mathcal{F}_a(\psi_k, A_k, L_{k + 1})$ \label{labelacf}
\STATE $\bm{y}_{k + 1} = \mathcal{F}_y(\bm{x}_k, \psi_k, A_k, a_{k + 1})$ \label{labelacg}
\STATE $ A_{k + 1} = A_{k} + a_{k + 1}$ \label{labelach}
\STATE $\bm{x}_{k + 1} = \displaystyle \argmin_{\bm{x} \in \mathbb{R}^n}u_{k + 1}(\bm{x})$ \label{labelaci}
\STATE $ \psi_{k + 1}(\bm{x}) = \psi_k(\bm{x}) + a_{k + 1} w_{k + 1}(\bm{x})$ \label{labelacj}
\ENDFOR
\end{algorithmic}
\end{pattern}

\subsection{Upper bounds}

Interestingly, all of the above results do not rely on a specific form of the local upper bound $u_{k + 1}$, as long as assumption \eqref{labelaba} holds for all $k \geq 0$. The enforced descent condition~\eqref{labelaba} can be \emph{equivalently} expressed in terms of composite objective values (by adding $\Psi(\bm{\bm{z}}_{k + 1})$ to both sides) as
\begin{equation} \label{labelack}
F(\bm{z}_{k + 1}) \leq Q_{L_{k + 1}, \bm{y}_{k + 1}}(\bm{z}_{k + 1}) + \Psi(\bm{z}_{k + 1}) , \quad k \geq 0.
\end{equation}
We want our algorithm to converge as fast as possible while maintaining the monotonicity property, expressed as
\begin{equation} \label{labelacl}
F(\bm{x}_{k + 1}) \leq F(\bm{x}_k) , \quad k \geq 0.
\end{equation}
Then, without further knowledge of the objective function, \eqref{labelack} and \eqref{labelacl} suggest a simple expression of the upper bound for all $\bm{x} \in \mathbb{R}^n$ and $k \geq 0$ in the form of
\begin{equation} \label{labelacm}
u_{k + 1}(\bm{x}) = \min\{ Q_{L_{k + 1}, \bm{y}_{k + 1}}(\bm{x}) + \Psi(\bm{x}), F(\bm{x}_k) + \sigma_{\{\bm{x}_k\}}(\bm{x}) \},
\end{equation}
where $\sigma_{X}$ is the indicator function~\cite{refabh} of set $X$, given by
\begin{equation}
\sigma_{X} (\bm{x}) = \left\{ \begin{array}{rl}
0 , &\bm{x} \in X, \\
+\infty, &\mbox{otherwise. }
\end{array} \right.
\end{equation}

\subsection{Towards an algorithm}

With the lower and upper bounds fully defined, we select functions $\mathcal{S}$, $\mathcal{F}_a$, and $\mathcal{F}_y$ so as to preserve the Lyapunov property of the gap sequence~\eqref{labelaby}.

\begin{theorem} \label{labelacn}
If at iteration $k \geq 0$ we have
\begin{equation} \nonumber
F(\bm{x}_{k + 1}) \leq F(\bm{z}_{k + 1}) \leq Q_{L_{k + 1}, \bm{y}_{k + 1}}(\bm{z}_{k + 1}) + \Psi(\bm{z}_{k + 1}),
\end{equation}
then
\begin{equation} \nonumber
\Delta_{k + 1} +\mathcal{A}_{k + 1} + \mathcal{B}_{k + 1} \leq \Delta_k,
\end{equation}
where subexpressions $\mathcal{A}_{k + 1}$, $\mathcal{B}_{k + 1}$, $\bm{g}_{k + 1}$, $\bm{s}_{k + 1}$, and $Y_{k + 1}$ are, respectively, defined as
\begin{equation} \nonumber
\begin{aligned}
\mathcal{A}_{k + 1} & \operatorname{\ \overset{def}{=} \ } \frac{1}{2} \left( \frac{A_{k + 1}}{L_{k + 1} + \mu_{\Psi}} - \frac{a_{k + 1}^2}{\gamma_{k + 1}} \right) \| \bm{g}_{k + 1} \|_2^2, \\
\mathcal{B}_{k + 1} & \operatorname{\ \overset{def}{=} \ } \frac{1}{ \gamma_{k + 1}} \langle \bm{g}_{k + 1} - \frac{\mu}{2 Y_{k + 1}} \bm{s}_{k + 1}, \bm{s}_{k + 1} \rangle, \\
\bm{g}_{k + 1} & \operatorname{\ \overset{def}{=} \ } (L_{k + 1} + \mu_{\Psi})(\bm{y}_{k + 1} - \bm{z}_{k + 1}), \\
\bm{s}_{k + 1} & \operatorname{\ \overset{def}{=} \ } A_k \gamma_{k + 1} \bm{x}_k + a_{k + 1} \gamma_k \bm{v}_k - Y_{k + 1} \bm{y}_{k + 1} , \\
Y_{k + 1} & \operatorname{\ \overset{def}{=} \ } A_k \gamma_{k + 1} + a_{k + 1} \gamma_k.
\end{aligned}
\end{equation}
\end{theorem}
\begin{proof}
We assume $k \geq 0$ throughout this proof. The descent condition assumption implies that the lower bound property in Lemma~\ref{labelaaz} holds. We define
\begin{equation} \label{labelaco}
\bm{G}_{k + 1} = \bm{g}_{k + 1} - \mu \bm{y}_{k + 1} = (L_{k + 1} - \mu_f) \bm{y}_{k + 1} - (L_{k + 1} + \mu_{\Psi}) \bm{z}_{k + 1},
\end{equation}
which simplifies the following subexpression of $w_{k + 1}$ in Lemma~\ref{labelaaz}:
\begin{equation} \label{labelacp}
\begin{gathered}
(L_{k + 1} + \mu_{\Psi}) \langle \bm{y}_{k + 1} - \bm{z}_{k + 1}, \bm{x} - \bm{y}_{k + 1} \rangle + \frac{\mu}{2}\| \bm{x} - \bm{y}_{k + 1} \|_2^2 \\
= \langle \bm{G}_{k + 1}, \bm{x} - \bm{y}_{k + 1} \rangle + \frac{\mu}{2} \| \bm{x} \|_2^2 - \frac{\mu}{2} \| \bm{y}_{k + 1} \|_2^2.
\end{gathered}
\end{equation}
The tightness of lower bound $w_{k + 1}$ in Lemma~\ref{labelaaz} is given by the residual $R_{k + 1}$ as
\begin{equation} \label{labelacq}
R_{k + 1}(\bm{x}) \operatorname{\ \overset{def}{=} \ } F(\bm{x}) - w_{k + 1}(\bm{x}) , \quad \bm{x} \in \mathbb{R}^n.
\end{equation}
From Lemma~\ref{labelaaz} it follows that
\begin{equation} \label{labelacr}
A_k R(\bm{x}_k) + a_{k + 1} R(\bm{x}^*) \geq 0.
\end{equation}
By expanding terms in \eqref{labelacr} using the identity in \eqref{labelacp} we obtain
\begin{equation} \label{labelacs}
A_{k} (F(\bm{x}_{k}) - F(\bm{x}^*)) - A_{k + 1} (F(\bm{z}_{k + 1}) - F(\bm{x}^*)) \geq \mathcal{C}_{k + 1},
\end{equation}
where
\begin{equation} \label{labelact}
\begin{gathered}
\mathcal{C}_{k + 1} \operatorname{\ \overset{def}{=} \ } \frac{A_{k + 1}}{L_{k + 1} + \mu_{\Psi}} \| \bm{g}_{k + 1} \|_2^2 + \langle \bm{G}_{k + 1}, A_k \bm{x}_k + a_{k + 1} \bm{x}^* - A_{k + 1} \bm{y}_{k + 1} \rangle \\
+ \frac{A_k \mu}{2} \| \bm{x}_k \|_2^2 + \frac{a_{k + 1} \mu}{2} \| \bm{x}^* \|_2^2 - \frac{A_{k + 1} \mu}{2} \| \bm{y}_{k + 1} \|_2^2.
\end{gathered}
\end{equation}
The vertex update in \eqref{labelabq} implies that
\begin{align}
a_{k + 1} \bm{G}_{k + 1} &= \gamma_k \bm{v}_k - \gamma_{k + 1} \bm{v}_{k + 1} , \label{labelacu} \\
a_{k + 1}^2 \| \bm{g}_{k + 1} \|_2^2 &= \| \gamma_k \bm{v}_k - \gamma_{k + 1} \bm{v}_{k + 1} + a_{k + 1} \mu \bm{y}_{k + 1}\|_2^2 . \label{labelacv}
\end{align}
Using \eqref{labelacu} and \eqref{labelacv} in \eqref{labelact} and rearranging terms yields
\begin{equation} \label{labelacw}
\mathcal{C}_{k + 1} = \mathcal{A}_{k + 1} + V_{k + 1} + \frac{1}{\gamma_{k + 1}} \langle \bm{G}_{k + 1}, \bm{s}_{k + 1} \rangle + \frac{\mu}{2 \gamma_{k + 1}} S_{k + 1},
\end{equation}
with $S_{k + 1}$ and $V_{k + 1}$ defined as
\begin{align}
S_{k + 1} &\ \operatorname{\ \overset{def}{=} \ } \ A_k \gamma_{k + 1} \| \bm{x}_k \|_2^2 + a_{k + 1} \gamma_k \| \bm{v}_k \|_2^2 - Y_{k + 1} \| \bm{y}_{k + 1} \|_2^2 , \\
V_{k + 1} &\ \operatorname{\ \overset{def}{=} \ } \ \frac{\gamma_{k + 1}}{2} \| \bm{v}_{k + 1} \|_2^2 - \frac{\gamma_k}{2} \| \bm{v}_k \|_2^2 + \langle \bm{G}_{k + 1}, a_{k + 1} \bm{x}^* \rangle + \frac{a_{k + 1} \mu}{2} \| \bm{x}^* \|_2^2 . \label{labelacx}
\end{align}
By further applying \eqref{labelabp} and \eqref{labelacu}, $S_{k + 1}$ and $V_{k + 1}$, respectively, become
\begin{align}
S_{k + 1} & = \left\langle \frac{1}{Y_{k + 1}} \bm{s}_{k + 1} + 2 \bm{y}_{k + 1}, \bm{s}_{k + 1} \right\rangle + \frac{a_{k + 1} A_k \gamma_k \gamma_{k + 1}}{ A_k \gamma_{k + 1} + a_{k + 1} \gamma_k } \| \bm{x}_k - \bm{v}_k \|_2^2 , \label{labelacy} \\
V_{k + 1} & = \frac{\gamma_{k + 1}}{2} \| \bm{v}_{k + 1} - \bm{x}^* \|_2^2 - \frac{\gamma_k}{2} \| \bm{v}_k - \bm{x}^* \|_2^2 , \label{labelacz}.
\end{align}
Putting together \eqref{labelaco}, \eqref{labelacs}, \eqref{labelacw}, \eqref{labelacy}, and \eqref{labelacz} with $F(\bm{x}_{k + 1}) \leq F(\bm{z}_{k + 1})$ and the fact that the second term in the right-hand side of \eqref{labelacy} is always non-negative, we obtain the desired result. More details can be found in proofs of Theorems 3, 4, and 5 in \cite{refabi}.
\end{proof}

Theorem~\ref{labelacn} provides a simple sufficient condition for the monotonicity of the gap sequence.
\begin{proposition} \label{labelada}
The monotonic decrease of the gap sequence in \eqref{labelaby} is guaranteed regardless of the algorithmic state if, for all $k \geq 0$, the following hold:
\begin{gather}
\bm{y}_{k + 1} = \frac{1}{Y_{k + 1}}(A_k \gamma_{k + 1} \bm{x}_k + a_{k + 1} \gamma_k \bm{v}_k) , \label{labeladb} \\
a_{k + 1} \leq \mathcal{E}(\gamma_k, A_k, L_{k + 1}) , \label{inlabeladf}
\end{gather}
where
\begin{equation}
\mathcal{E}(\gamma_k, A_k, L_{k + 1}) \operatorname{\ \overset{def}{=} \ } \frac{\gamma_k + A_k \mu}{2 (L_{k + 1} - \mu_f)} \left(1 + \sqrt{1 + \frac{4 (L_{k + 1} - \mu_f) A_k \gamma_k}{(\gamma_k + A_k \mu)^2}} \right).
\end{equation}
\end{proposition}
\begin{proof}
Herein, we consider all $k \geq 0$. According to Theorem~\ref{labelacn}, an obvious sufficient condition for \eqref{labelaby} is $\mathcal{A}_{k + 1} \geq 0$ combined with $\mathcal{B}_{k + 1} \geq 0$. Due to the black-box nature of the oracle functions and the assumption of arbitrary algorithmic state, we allow for $\bm{g}_{k + 1}$ to be any vector in $\mathbb{R}^n$.
Quantity $\mathcal{B}_{k + 1}$ is a scalar product between two vectors, the one containing $\bm{g}_{k + 1}$ being arbitrary. Therefore, a sufficient condition for $\mathcal{B}_{k + 1} \geq 0$ is that $\bm{s}_{k + 1} = 0$, which we can always guarantee by setting \eqref{labeladb}.
Moreover, the value of $\| \bm{g}_{k + 1} \|_2^2$ could be very large, but is always non-negative. Therefore, a sufficient condition for $\mathcal{A}_{k + 1} \geq 0$ is that
\begin{equation} \label{labeladd}
(L_{k + 1} + \mu_{\Psi}) a_{k + 1}^2 \leq A_{k + 1} \gamma_{k + 1} .
\end{equation}
By expanding \eqref{labeladd} using \eqref{labelabh} and \eqref{labelabp} and rearranging terms, we obtain that $a_{k + 1}$ needs to satisfy the inequality
\begin{equation} \label{labelade}
(L_{k + 1} - \mu_f) a_{k + 1}^2 - (\gamma_k + A_k \mu) a_{k + 1} - A_k \gamma_k \leq 0.
\end{equation}
Because $L_{k + 1} > \mu_f$, the solutions of \eqref{labelade} lie between the two roots of the corresponding equation, only one of which is positive and given by $\mathcal{E}(\gamma_k, A_k, L_{k + 1})$.
\end{proof}
Proposition~\ref{labelada} allows us to select $\mathcal{F}_a$ and $\mathcal{F}_y$ as to preserve the Lyapunov property in \eqref{labelaby}. For $\mathcal{F}_y$, the obvious choice is \eqref{labeladb}. Out of the multitude of potential $\mathcal{F}_a$, we choose the one that yields the largest convergence guarantees $A_{k + 1}$, given by
\begin{equation} \label{labeladf}
\mathcal{F}_a(\psi_k, A_k, L_{k + 1}) = \mathcal{E}(\gamma_k, A_k, L_{k + 1}), \quad k \geq 0,
\end{equation}
which corresponds to equality in \eqref{inlabeladf}, namely
\begin{equation} \label{labeladg}
(L_{k + 1} + \mu_{\Psi}) a^2_{k + 1} = A_{k + 1} \gamma_{k + 1}, \quad k \geq 0.
\end{equation}

For determining the LCE, we select the backtracking line-search method $\mathcal{S}_A$ employed by AMGS ~\cite{refaak} and the original ACGM~\cite{refaaz}. The search parameters comprise the LCE increase rate $r_u > 1$ and the LCE decrease rate $0 < r_d < 1$. The search terminates when the line-search stopping criterion (LSSC) in \eqref{labelaba} is satisfied.

\subsection{Putting it all together}

We have thus determined a search strategy $\mathcal{S}_A$, the initial estimate function $\psi_0$ in~\eqref{labelabj}, the upper bounds $u_{k + 1}$ in~\eqref{labelacm}, the lower bounds $w_{k + 1}$ in Lemma~\ref{labelaaz}, function $\mathcal{F}_a$ in~\eqref{labeladf}, and function $\mathcal{F}_y$ in~\eqref{labeladb}. Substituting these expressions in the design model outlined in Pattern~\ref{labelacc}, we can write down a generalization of ACGM in estimate sequence form, as listed in Algorithm~\ref{labeladi}.

Non-monotone generalized ACGM can be obtained by enforcing $\bm{x}_{k + 1} = \bm{z}_{k + 1}$ for all $k \geq 0$, accomplished by replacing \cmdaah{labeladi}{labelads} with
\begin{equation} \label{labeladh}
\bm{x}_{k + 1} := \cmdaac{z}.
\end{equation}

\begin{algorithm}[h!]
\caption{Generalized monotone ACGM in estimate sequence form \newline
\textbf{ACGM}($\bm{x_0}$, $L_0$, $\mu_f$, $\mu_{\Psi}$, $A_0$, $\gamma_0$, $r_u$, $r_d$, $K$)}
\label{labeladi}
\begin{algorithmic}[1]
\STATE $\bm{v}_0 = \bm{x}_0, \ \mu = \mu_f + \mu_{\Psi}$
\FOR {$k = 0,\ldots{},K-1$}
\STATE $\cmdaab{L} := r_d L_k$ \label{labeladi_start_search}
\LOOP
\STATE $\cmdaab{a} := \frac{\gamma_k + A_k \mu}{2 (\cmdaab{L} - \mu_f)} \left(1 + \sqrt{1 + \frac{4 (\cmdaab{L} - \mu_f) A_k \gamma_k}{(\gamma_k + A_k \mu)^2}} \right)$ \label{labeladk}
\STATE $\cmdaab{A} := A_k + \cmdaab{a}$ \label{labeladl}
\STATE $\cmdaab{\gamma} := \gamma_k + \cmdaab{a} \mu$ \label{labeladm}
\STATE $\cmdaac{y} := \frac{1}{A_k \cmdaab{\gamma} +
\cmdaab{a} \gamma_k} (A_k \cmdaab{\gamma} \bm{x}_k + \cmdaab{a} \gamma_k \bm{v}_k)$ \label{labeladn}
\STATE $\cmdaac{z} :=\mbox{prox}_{\frac{1}{\cmdaab{L}} \Psi}\left(\cmdaac{y}
- \frac{1}{\cmdaab{L}}\nabla f(\cmdaac{y}) \right)$ \label{labelado}
\IF {$f(\cmdaac{z}) \leq Q_{\cmdaab{L}, \cmdaac{y}}(\cmdaac{z})$} \label{labeladp}
\STATE Break from loop \label{labeladq}
\ELSE
\STATE $\cmdaab{L} := r_u \cmdaab{L} $
\ENDIF
\ENDLOOP \label{labeladr}
\STATE $\bm{x}_{k + 1} := \argmin \{ F(\cmdaac{z}), F(\bm{x}_k) \}$ \label{labelads}
\STATE $\bm{v}_{k + 1} := \frac{1}{\cmdaab{\gamma}} (\gamma_k \bm{v}_k + \cmdaab{a} (\cmdaab{L} + \mu_{\Psi}) \cmdaac{z} - \cmdaab{a}(\cmdaab{L} - \mu_{f}) \cmdaac{y})$ \label{labeladt}
\STATE $L_{k + 1} := \cmdaab{L}, \ A_{k + 1} := \cmdaab{A}, \ \gamma_{k + 1} := \cmdaab{\gamma}$

\ENDFOR
\STATE \textbf{return} $x_K$
\end{algorithmic}
\end{algorithm}

\section{Complexity analysis}

\subsection{Worst-case convergence guarantees}

Algorithm~\ref{labeladi} maintains the convergence guarantee in \eqref{labelaam} explicitly at run-time as state variable $A_k$. Moreover, if sufficient knowledge of the problem is available, it is possible to formulate a worst-case convergence guarantee \emph{before} running the algorithm.

For our analysis, we will need to define a number of curvature-related quantities, namely the local inverse condition number $q_{k + 1}$ for all $k \geq 0$, the worst-case LCE $L_u$, and the worst-case inverse condition number $q_u$, respectively, given by
\begin{gather}
q_{k + 1} \operatorname{\ \overset{def}{=} \ } \frac{\mu}{L_{k + 1} + \mu_{\Psi}}, \label{labeladu} \\
L_u \operatorname{\ \overset{def}{=} \ } \max\{r_u L_f, r_d L_0\}, \\
q_u \operatorname{\ \overset{def}{=} \ } \frac{\mu}{L_u + \mu_{\Psi}}.
\end{gather}

The worst-case convergence rate for generalized ACGM is stated in the following theorem.
\begin{theorem} \label{labeladv}
If $\gamma_0 \geq A_0 \mu$, the generalized ACGM algorithm generates a sequence $\{\bm{x}_k\}_{k \geq 1}$ that satisfies
\begin{equation} \nonumber
F(\bm{x}_k) - F(\bm{x}^*) \leq \min\left\{\frac{4}{(k + 1)^2}, (1 - \sqrt{q_u})^{k - 1}\right\}
(L_u - \mu_f) \bar{\Delta}_0 , \quad k \geq 1,
\end{equation}
where
\begin{equation} \nonumber
\bar{\Delta}_0 \operatorname{\ \overset{def}{=} \ } \frac{\Delta_0}{\gamma_0} = \frac{A_0}{\gamma_0} (F(\bm{x}_0) - F(\bm{x}^*)) + \frac{1}{2} \| \bm{x}_0 - \bm{x}^* \|_2^2.
\end{equation}
\end{theorem}
\begin{proof}
The non-negativity of the weights~\eqref{labelabi} implies that $\gamma_k \geq \gamma_0$ for all $k \geq 0$. Combined with \eqref{labeladg} and $A_k \mu \geq 0$, we have for all $k \geq 0$ that
\begin{equation}
\begin{gathered}
A_{k + 1} = A_k + a_{k + 1} \geq A_k + \frac{\gamma_0}{2 (L_{k + 1} - \mu_f)}
+ \sqrt{\frac{\gamma_0^2}{4 (L_{k + 1} - \mu_f)^2} + \frac{A_k \gamma_0}{(L_{k + 1} - \mu_f)}}.
\end{gathered}
\end{equation}
As we can see from Algorithm~\ref{labeladi}, scaling $A_0$ and $\gamma_0$ by a fixed factor does not alter the behavior of generalized ACGM. Additionally, $\gamma_0$ is guaranteed to be non-zero. To simplify calculations, we introduce the normalized convergence guarantees $\bar{A}_k \operatorname{\ \overset{def}{=} \ } A_k / \gamma_0$ for all $k \geq 0$.

Regardless of the outcome of individual line-search calls, the growth of the normalized accumulated weights obeys
\begin{equation}
\bar{A}_{k + 1} \geq \bar{A}_k +
\frac{1}{2 (L_u - \mu_f)} +
\sqrt{\frac{1}{4 (L_u - \mu_f)^2} + \frac{\bar{A}_k}{(L_u - \mu_f)}}, \quad k \geq 0.
\end{equation}
Taking into account that $A_0 \geq 0$, we obtain by induction that
\begin{equation} \label{labeladw}
\bar{A}_k \geq \frac{(k + 1)^2}{4 (L_u - \mu_f)}, \quad k \geq 1.
\end{equation}

From assumption $\gamma_0 \geq A_0 \mu$, \eqref{labelabs} implies $\gamma_k \geq A_k \mu$ for all $k \geq 0$. Hence
\begin{equation}
\frac{a_{k + 1}^2}{A_{k + 1}^2} \overset{\eqref{labeladg}}{=} \frac{\gamma_{k + 1}}{(L_{k + 1} + \mu_{\Psi}) A_{k + 1}} \geq \frac{\mu}{L_{k + 1} + \mu_{\Psi}} = q_{k + 1} \geq q_u, \quad k \geq 0.
\end{equation}
Since $A_0 \geq 0$, we have that $\bar{A}_1 \geq \frac{1}{L_u - \mu_f}$. By induction, it follows that
\begin{equation} \label{labeladx}
\bar{A}_k \geq \frac{1}{L_u - \mu_f} (1 - \sqrt{q_u})^{-(k - 1)} , \quad k \geq 1.
\end{equation}
Substituting \eqref{labeladw} and \eqref{labeladx} in \eqref{labelaam} completes the proof.
\end{proof}

Theorem~\ref{labeladv} shows that generalized ACGM has the best convergence guarantees among its class of algorithms (see~\cite{refaaz} for an in-depth discussion).

Note that the assumption $\gamma_0 \geq A_0 \mu$ always holds for non-strongly convex objectives and that ACGM is guaranteed converge for strongly convex objectives also when $\gamma_0 < A_0 \mu$. However, in the latter case, it is more difficult to obtain simple lower bounds on the convergence guarantees. We leave such an endeavor to future research.

\subsection{Wall-clock time units} \label{labelady}

So far, we have measured the theoretical performance of algorithms in terms of convergence guarantees (including the worst-case ones) indexed in iterations. This does not account for the complexity of individual iterations. In~\cite{refaaz}, we have introduced a new measure of complexity, the \emph{wall-clock time unit} (WTU), to compare optimization algorithms more reliably. We thus distinguish between two types of convergence guarantees. One is the previously used \emph{iteration convergence guarantee}, indexed in iterations and a new \emph{computational convergence guarantee}, indexed in WTU.

The WTU is a measure of running time in a \emph{shared memory parallel scenario}. The computing environment consists of a small number of parallel processing units (PPU). Each PPU may be a parallel machine itself. The number of parallel units is considered sufficient to compute any number of independent oracle functions simultaneously. The shared-memory system does not impose constraints on parallelization, namely, it is uniform memory access (UMA)~\cite{refabj} and it is large enough to store the arguments and results of oracle calls for as long as they are needed.

In order to compare algorithms based on a unified benchmark, in~\cite{refaaz} we have assumed that $f(\bm{x})$ and $\nabla f(\bm{x})$ require $1$ WTU each while all other operations are negligible and amount to $0$~WTU. In this work, we generalize the analysis. We attribute finite non-negative costs $t_f$, $t_g$, $t_{\Psi}$, and $t_p$ to $f(\bm{x})$, $\nabla f(\bm{x})$, $\Psi(\bm{x} )$, and $\mbox{prox}_{\tau \Psi}(\bm{x})$, respectively. However, since we are dealing with large-scale problems, we maintain the assumption that element-wise vector operations, including scalar-vector multiplications, vector additions, and inner products, have negligible complexity when compared to oracle functions and assign a cost of $0$ WTU to each. Synchronization of PPUs also incurs no cost. Consequently, when computed in isolation, an objective function value call costs $t_F = \max\{t_f, t_{\Psi}\}$, ascribable to separability, while a proximal gradient operation costs $t_T = t_g + t_p$, due to computational dependencies.

\subsection{Per-iteration complexity} \label{labeladz}

We measure this complexity in WTU on the shared memory system described in the previous subsection and consider a parallel implementation involving speculative execution~\cite{refabj}.

The advancement phase of a generalized ACGM iteration consists of one proximal gradient step (\cmdaah{labeladi}{labelado}). Hence, every iteration has a base cost of $t_T = t_g + t_p$. LSSC and the monotonicity condition (MC) in \cmdaah{labeladi}{labelads} can be evaluated in parallel with subsequent iterations. Both rely on the computation of $f(\cmdaac{z})$, which in the worst case requires $\lceil t_f / t_T \rceil$ dedicated PPUs. In addition, MC may need up to $\lceil t_{\Psi} / t_T \rceil$ PPUs.

Backtracks stall the algorithm in a way that cannot be alleviated by parallelization or intensity reduction. Therefore, it is desirable to make them a rare event. Assuming that the local curvature of $f$ varies around a fixed value, this would mean that $\log(r_u)$ should be significantly larger than $-log(r_d)$. With such a parameter choice, the algorithm can proceed from one iteration to another by speculating that backtracks do not occur at all. Let the current iteration be indexed by $k$. If the LSSC of iteration $k$ fails, then the algorithm discards all the state information pertaining to all iterations made after $k$, reverts to iteration $k$, and performs the necessary computation to correct the error. We consider that a misprediction incurs a detection cost $t_d$ and a correction cost $t_c$. LSSC requires the evaluation of $f(\cmdaac{z})$ and incurs a detection cost of $t_d = t_f$. A backtrack entails recomputing $\cmdaac{y}$, yielding an LSSC $t_c = t_T$ correction time.

Overshoots are assumed to occur even less often. Similarly, the algorithm proceeds speculating that MC always passes and defaults to~\eqref{labeladh}. Hence, MC has $t_d = t_F$, due to its dependency on $\Psi(\cmdaac{z})$, but once the algorithmic state of iteration $k$ has been restored, no additional oracle calls are needed, leading to $t_c = 0$. MC and LSSC can be fused into a single condition, giving rise to the scenarios outlined in Table~\ref{labelaea}. Note that if LSSC fails, MC is not evaluated.

\begin{table}[h!]
\centering
\small
\caption{Algorithm stall time in WTU based on the outcome of LSSC and MC}
\label{labelaea}
\begin{tabular}{c c c} \toprule
& MC passed & MC failed \\ \midrule
LSSC passed & $0$ & $\max\{t_f, t_{\Psi}\}$ \\
LSSC failed & $t_f + t_g + t_p$ & N / A \\ \bottomrule
\end{tabular}
\end{table}

For non-monotone generalized ACGM, each backtrack adds $t_f + t_T$ WTU to a base iteration cost of $t_T$. A comparison to other methods employing line-search is shown in Table~\ref{labelaeb}.

\begin{table}[h!]
\centering
\small
\caption{Per-iteration cost of FISTA, AMGS, and generalized ACGM in the non-monotone setting}
\label{labelaeb}
\begin{tabular}{l c c c} \toprule
& FISTA & AMGS & ACGM \\ \midrule
Base cost & $t_g + t_p$ & $2 t_g + 2 t_p$ & $t_g + t_p$ \\
LSSC $t_d$ & $t_f$ & $t_g$ & $t_f$ \\
LSSC $t_c$ & $t_p$ & $t_g + t_p$ & $t_g + t_p$ \\
Backtrack cost & $t_f + t_p$ & $2 t_g + t_p$ & $t_f + t_g + t_p$ \\ \bottomrule
\end{tabular}
\end{table}

\clearpage

\section{Extrapolated form}

\subsection{Monotonicity and extrapolation}

In the original ACGM~\cite{refaaz}, the auxiliary point can be obtained from two successive main iterates through extrapolation. Interestingly, this property is preserved for any step size. We show in the following how monotonicity alters this property and bring generalized monotone ACGM to a form in which the auxiliary point is an extrapolation of state variables. We begin with the following observation.
\begin{lemma} \label{labelaec}
In Algorithm~\ref{labeladi}, the estimate sequence vertices can be obtained from other state variables through extrapolation as
\begin{equation} \nonumber
\bm{v}_{k + 1} = \bm{x}_k + \frac{A_{k + 1}}{a_{k + 1}} (\bm{z}_{k + 1} - \bm{x}_k) , \quad k \geq 0.
\end{equation}
\end{lemma}
\begin{proof}
In this proof we consider all $k \geq 0$. The vertex update in \eqref{labelabq} can be rewritten using \eqref{labeladg} as
\begin{align}
a_{k + 1} \gamma_{k} \bm{v}_k &\overset{\eqref{labelabq}}{=} a_{k + 1} \gamma_{k + 1} \bm{v}_{k + 1} + (L_{k + 1} + \mu_{\Psi}) a^2_{k + 1} (\bm{y}_{k + 1} - \bm{z}_{k + 1}) + a_{k + 1}^2 \mu \bm{y}_{k + 1} \nonumber \\
&\overset{\eqref{labeladg}}{=} a_{k + 1} \gamma_{k + 1} \bm{v}_{k + 1} + A_{k + 1} \gamma_{k + 1} (\bm{y}_{k + 1} - \bm{z}_{k + 1}) + a_{k + 1}^2 \mu \bm{y}_{k + 1} . \label{labelaed}
\end{align}
The auxiliary point update in \eqref{labeladb} can also be written as
\begin{equation} \label{labelaee}
a_{k + 1} \gamma_{k} \bm{v}_k = (A_k \gamma_{k + 1} + a_{k + 1} \gamma_k) \bm{y}_{k + 1} - A_k \gamma_{k + 1} \bm{x}_k.
\end{equation}
From \eqref{labelabp} we have that
\begin{equation} \label{labelaef}
A_{k + 1} \gamma_{k + 1} - a_{k + 1}^2 \mu = A_k \gamma_{k + 1} + a_{k + 1} \gamma_k.
\end{equation}
Putting together \eqref{labelaed}, \eqref{labelaee}, and eliminating $\bm{y}_{k + 1}$ using \eqref{labelaef} we obtain
\begin{equation}
-A_k \gamma_{k + 1} \bm{x}_k = a_{k + 1} \gamma_{k + 1} \bm{v}_{k + 1} - A \gamma_{k + 1} \bm{z}_{k + 1}.
\end{equation}
\end{proof}
Lemma~\ref{labelaec} enables us to write down the corresponding expression for the auxiliary point.
\begin{proposition} \label{labelaeg}
In Algorithm~\ref{labeladi}, the auxiliary points $\bm{y}_{k + 1}$ obey the following extrapolation rule:
\begin{equation} \nonumber
\bm{y}_{k + 1} = \bm{x}_k + \beta_k (\bm{z}_k - \bm{x}_{k - 1}) , \quad k \geq 1,
\end{equation}
where
\begin{align}
\beta_k &= \left(\frac{A_k}{a_k} - \bm{1}_{\{\bm{z}_k\}}(\bm{x}_k) \right) \omega_k , \label{labelaeh}\\
\omega_k &\operatorname{\ \overset{def}{=} \ } \frac{ a_{k + 1} \gamma_k }{A_k \gamma_{k + 1} + a_{k + 1} \gamma_k },
\end{align}
and
$\bm{1}_X$ denotes the membership function of set $X$, namely
\begin{equation} \nonumber
\bm{1}_X(x) = \left\{ \begin{array}{rl} 1, & x \in X \\ 0, & x \notin X \end{array} \right..
\end{equation}
\end{proposition}
\begin{proof}
The results herein apply to all $k \geq 1$. Lemma~\ref{labelaec} applied for $\bm{v}_k$, combined with the auxiliary point update in \eqref{labeladb}, leads to
\begin{equation} \label{labelaei}
\begin{gathered}
\bm{y}_{k + 1} = (1 - \omega_k) \bm{x}_k + \frac{A_k}{a_k} \omega_k \bm{z}_k + \left(1 - \frac{A_k}{a_k} \right) \omega_k \bm{x}_{k - 1}.
\end{gathered}
\end{equation}
Depending on the outcome of the update in \cmdaah{labeladi}{labelads}, we distinguish two situations.

If MC passes at iteration $k - 1$ ($F(\bm{z}_k) \leq F(\bm{x}_{k - 1})$), we set $\bm{x}_k = \bm{z}_k$, hence
\begin{equation} \label{labelaej}
\bm{y}_{k + 1} = \bm{x}_k + \left(\frac{A_k}{a_k} - 1 \right) \omega_k (\bm{z}_k - \bm{x}_{k - 1}).
\end{equation}

If the algorithm overshoots ($F(\bm{z}_k) > F(\bm{x}_{k - 1})$) then, by monotonicity, we impose $\bm{x}_k = \bm{x}_{k - 1}$, which leads to
\begin{equation} \label{labelaek}
\bm{y}_{k + 1} = \bm{x}_k + \frac{A_k}{a_k} \omega_k (\bm{z}_k - \bm{x}_{k - 1}).
\end{equation}
Combining \eqref{labelaej} and \eqref{labelaek} completes the proof.
\end{proof}

Until this point we have assumed that the first iteration $k = 0$ \emph{does not} use the auxiliary point extrapolation rule in Proposition~\ref{labelaeg}. To write our generalized ACGM in a form similar to monotone FISTA (MFISTA~\cite{refaba}) and the monotone version of FISTA-CP~\cite{refaah}, we define the vertex extrapolation factor in Lemma~\ref{labelaec} as
\begin{equation} \label{labelael}
t_k \operatorname{\ \overset{def}{=} \ } \frac{A_k}{a_k} , \quad k \geq 1.
\end{equation}
\begin{proposition} \label{labelaem}
For $k \geq 1$, the vertex extrapolation factor can be determined using a recursion rule that does not depend on weights $a_k$ and $A_k$, given by
\begin{equation} \label{labelaen}
t_{k + 1}^2 + t_{k + 1} (q_k t_k^2 - 1) - \frac{L_{k + 1} + \mu_{\Psi}}{L_k + \mu_{\Psi}}t_k^2 = 0.
\end{equation}
Subexpression $\omega_k$ and auxiliary point extrapolation factor $\beta_k$ can also be written for all $k \geq 1$ as
\begin{gather}
\omega_k = \frac{1 - q_{k + 1} t_{k + 1}}{(1 - q_{k + 1}) t_{k + 1}}, \label{labelaeo} \\
\beta_k = (t_k - \bm{1}_{\{\bm{z}_k\}}(\bm{x}_k)) \omega_k. \label{labelaep}
\end{gather}
\end{proposition}
\begin{proof}
In this proof we take all $k \geq 1$.
From \eqref{labeladu} it follows that
\begin{equation} \label{labelaeq}
q_k = \frac{L_{k + 1} + \mu_{\Psi}}{L_k + \mu_{\Psi}} q_{k + 1}.
\end{equation}
Moreover,
\begin{align}
\frac{\gamma_k}{\gamma_{k + 1}} &\overset{\eqref{labelabp}}{=} \frac{A_{k + 1}(\gamma_{k + 1} - a_{k + 1} \mu)}{A_{k + 1} \gamma_{k + 1}} \overset{\eqref{labeladg}}{=} 1 - \frac{A_{k + 1} a_{k + 1} \mu}{(L_{k + 1} + \mu_{\Psi}) a_{k + 1}^2} \overset{\eqref{labelael}}{=} 1 - q_{k + 1} t_{k + 1} , \label{labelaer} \\
\gamma_{k} t_k^2 &\overset{\eqref{labelael}}{=} \frac{(L_k + \mu_{\Psi})\gamma_{k} A_k^2}{(L_k + \mu_{\Psi}) a_k^2} \overset{\eqref{labeladg}}{=} (L_k + \mu_{\Psi}) A_{k} , \label{labelaes} \\
\gamma_{k} t_k &\overset{\eqref{labelael}}{=} \frac{\gamma_k A_k}{a_k} \overset{\eqref{labeladg}}{=} (L_k + \mu_{\Psi}) a_{k} . \label{labelaet}
\end{align}
Putting together \eqref{labelaeq}, \eqref{labelaer}, \eqref{labelaes}, and \eqref{labelaet} we obtain
\begin{equation} \nonumber
\begin{gathered}
\gamma_{k + 1} \left( t_{k + 1}^2 + t_{k + 1} (q_k t_k^2 - 1) - \frac{L_{k + 1} + \mu_{\Psi}}{L_k + \mu_{\Psi}}t_k^2 \right) \\
\overset{\eqref{labelaeq}}{=} \gamma_{k + 1} t_{k + 1}^2 + \frac{L_{k + 1} + \mu_{\Psi}}{L_k + \mu_{\Psi}} \gamma_{k + 1} q_{k + 1} t_{k + 1} t_k^2 - \gamma_{k + 1} t_{k + 1} - \gamma_{k + 1} \frac{L_{k + 1} + \mu_{\Psi}}{L_k + \mu_{\Psi}} t_k^2 \\
\overset{\eqref{labelaer}}{=} (L_{k + 1} + \mu_{\Psi}) A_{k + 1} - \frac{L_{k + 1} + \mu_{\Psi}}{L_k + \mu_{\Psi}} \gamma_k t_k^2 - (L_{k + 1} + \mu_{\Psi}) a_{k + 1} \\
= (L_{k + 1} + \mu_{\Psi}) A_{k + 1} - (L_{k + 1} + \mu_{\Psi}) A_{k} - (L_{k + 1} + \mu_{\Psi}) a_{k + 1} = 0.
\end{gathered}
\end{equation}
Since $\gamma_{k + 1} > 0$, \eqref{labelaen} holds. Subexpression $\omega_k$ is similarly obtained as
\begin{equation} \nonumber
\omega_k \overset{\eqref{labelaef}}{=} \frac{1 - \frac{a_{k + 1} \mu}{\gamma_{k + 1}}}{\frac{A_{k + 1}}{a_{k + 1}} - \frac{a_{k + 1} \mu}{\gamma_{k + 1}} } \overset{\eqref{labeladg}}{=} \frac{1 - \frac{\mu A_{k + 1}}{(L_{k + 1} + \mu_{\Psi}) a_{k + 1}}}{\left(1 - \frac{\mu}{L_{k + 1} + \mu_{\Psi}} \right) \frac{A_{k + 1}}{a_{k + 1}}} \overset{\eqref{labelael}}{=} \frac{1 - q_{k + 1} t_{k + 1}}{(1 - q_{k + 1}) t_{k + 1}}.
\end{equation}
Substituting \eqref{labelael} in \eqref{labelaeh} gives \eqref{labelaep}.
\end{proof}

For simplicity, we wish to extend the update rules from Propositions \ref{labelaeg} and \ref{labelaem} to the first iteration $k = 0$. The missing parameters follow naturally from this extension. First, $t_0$ can be obtained by setting $k = 0$ in \eqref{labelaen} as
\begin{equation} \label{labelaeu}
t_0 = \sqrt{\frac{t_1^2 - t_1}{\frac{L_1 + \mu_{\Psi}}{L_0 + \mu_{\Psi}} - t_1 q_0}}
\overset{\eqref{labelael}}{=} \sqrt{\frac{A_1 - a_1}{\frac{(L_1 + \mu_{\Psi}) a_1^2}{(L_0 + \mu_{\Psi}) A_1} - a_1 q_0}} \overset{\eqref{labeladg}}{=} \sqrt{\frac{(L_0 + \mu_{\Psi}) A_0}{\gamma_0}}.
\end{equation}
Next, we introduce a ``phantom iteration'' $k = -1$ with the main iterate as the only state parameter. We set $x_{-1} \operatorname{\overset{def}{=}} x_0$ so that \emph{any} value of $\beta_0$ will satisfy Proposition~\ref{labelaeg}. For brevity, we choose to obtain $\beta_0$ from expression \eqref{labelaep} with $k = 0$.

Thus, with the initialization in \eqref{labelaeu} and the recursion in \eqref{labelaen}, we have completely defined the vertex extrapolation factor sequence $\{t_k\}_{k \geq 0}$, and derived from it the auxiliary extrapolation factor expression in \eqref{labelaep}. Now, we no longer need to maintain weight sequences $\{a_k\}_{k \geq 1}$ and $\{A_k\}_{k \geq 0}$. We simplify our generalized ACGM further by noting that, to produce the auxiliary point, the extrapolation rule in Proposition~\ref{labelaeg} depends on three vector parameters. However, it is not necessary to store both $\bm{z}_k$ and $\bm{x}_{k - 1}$ across iterations. To address applications where memory is limited, we only maintain the difference term $\bm{d}_k$, given by
\begin{equation}
\bm{d}_k = (t_k - \bm{1}_{\{\bm{z}_k\}}(\bm{x}_k)) (\bm{z}_k - \bm{x}_{k - 1}), \quad k \geq 0. \label{labelaev}
\end{equation}

The extrapolation rule in Proposition~\ref{labelaeg} becomes
\begin{equation}
\bm{y}_{k + 1} = \bm{x}_k + \omega_k \bm{d}_k, \quad k \geq 0.
\end{equation}
Note that subexpression $\omega_k$ contains only recent information whereas $\bm{d}_k$ needs only to access the state of the preceding iterations.

The above modifications yield a form of generalized ACGM based on extrapolation, which we list in Algorithm~\ref{labelaew}. To obtain a non-monotone algorithm, it suffices to replace \cmdaah{labelaew}{labelafe} with \eqref{labeladh}.

We stress that while Algorithms \ref{labeladi} and \ref{labelaew} carry out different computations, they are mathematically equivalent with respect to the main iterate sequence $\{\bm{x}_k\}_{k \geq 0}$. The oracle calls and their dependencies in Algorithm~\ref{labelaew} are also identical to those in Algorithm~\ref{labeladi}. Therefore the per-iteration complexity is the same.

\begin{algorithm}[h!]
\caption{Generalized monotone ACGM in extrapolated form \newline
\textbf{ACGM}($\bm{x_0}$, $L_0$, $\mu_f$, $\mu_{\Psi}$, $A_0$, $\gamma_0$, $r_u$, $r_d$, $K$)}
\label{labelaew}
\begin{algorithmic}[1]
\STATE $\bm{x}_{-1} = \bm{x}_0, \ \bm{d}_0 = \bm{0}$
\STATE $\mu = \mu_f + \mu_{\Psi}, \ t_0 = \sqrt{\frac{(L_0 + \mu_{\Psi}) A_0}{\gamma_0}}, \ q_0 = \frac{\mu}{L_0 + \mu_{\Psi}}$\\[1mm]
\FOR {$k = 0, ..., K - 1$}
\STATE $\cmdaab{L} := r_d L_k$ \label{labelaex}
\LOOP
\STATE $\cmdaab{q} := \frac{\mu}{\cmdaab{L} + \mu_{\Psi}}$ \label{labelaey}
\STATE $\cmdaab{t} := \frac{1}{2} \left(1 - q_k t_k^2 + \sqrt{(1 - q_k t_k^2)^2 + 4 \frac{\cmdaab{L} + \mu_{\Psi}}{L_k + \mu_{\Psi}} t_k^2} \right)$ \label{labelaez}
\STATE $\cmdaac{y} := \bm{x}_k + \frac{1 - \cmdaab{q} \cmdaab{t}}{(1 - \cmdaab{q}) \cmdaab{t}} \bm{d}_k$ \label{labelafa}
\STATE $\cmdaac{z} :=\mbox{prox}_{\frac{1}{\cmdaab{L}} \Psi}\left(\cmdaac{y} - \frac{1}{\cmdaab{L}}\nabla f(\cmdaac{y}) \right)$ \label{labelaew_update_z}

\IF {$f(\cmdaac{z}) \leq Q_{\cmdaab{L}, \cmdaac{y}}(\cmdaac{z})$} \label{labelafc}
\STATE Break from loop \label{labelaew_found_search}
\ELSE
\STATE $\cmdaab{L} := r_u \cmdaab{L} $
\ENDIF
\ENDLOOP
\STATE $\bm{x}_{k + 1} := \argmin \{ F(\cmdaac{z}), F(\bm{x}_k) \}$ \label{labelafe}
\STATE $\bm{d}_{k + 1} := (\cmdaab{t} - \bm{1}_{\{\cmdaac{z}\}}(\bm{x}_{k + 1})) (\cmdaac{z} - \bm{x}_k)$
\STATE $L_{k + 1} := \cmdaab{L}, \ q_{k + 1} := \cmdaab{q}, \ t_{k + 1} := \cmdaab{t}$
\ENDFOR
\STATE \textbf{return} $x_K$
\end{algorithmic}
\end{algorithm}

\subsection{Retrieving the convergence guarantee}

For Algorithm~\ref{labeladi}, the convergence guarantee in \eqref{labelaam} is obtained directly from the state variable $A_k$. For Algorithm~\ref{labelaew}, we need the following result.

\begin{lemma} \label{labelaff}
The vertex extrapolation factor expression in \eqref{labelaeu} generalizes to arbitrary $k \geq 0$ as
\begin{equation} \nonumber
t_k = \sqrt{\frac{(L_k + \mu_{\Psi}) A_k}{\gamma_k}}.
\end{equation}
\end{lemma}
\begin{proof}
For $k = 0$, \eqref{labelaeu} implies that Lemma~\ref{labelaff} holds.
For $k \geq 1$ we have that
\begin{equation} \nonumber
t_{k + 1} = \sqrt{\frac{(L_{k + 1} + \mu_{\Psi}) A_{k + 1}^2}{(L_{k + 1} + \mu_{\Psi}) a_{k + 1}^2}} \overset{\eqref{labeladg}}{=} \sqrt{\frac{(L_{k + 1} + \mu_{\Psi}) A_{k + 1}}{\gamma_{k + 1}}}.
\end{equation}
\end{proof}

From Lemma~\ref{labelaff}, we distinguish two scenarios.

If $\gamma_0 \neq A_0 \mu$, from \eqref{labelabs}, \eqref{labeladu}, and Lemma~\ref{labelaff} it follows that the convergence guarantee can be derived directly from the state parameters without alterations to Algorithm~\ref{labelaew} as
\begin{equation}
A_k = \frac{(\gamma_0 - A_0 \mu) t_k^2}{(L_k + \mu_{\Psi}) (1 - q_k t_k^2)}, \quad k \geq 1.
\end{equation}

\subsection{Border-case}

However, if $\gamma_0 = A_0 \mu$, then \eqref{labelabs}, \eqref{labeladu}, and Lemma~\ref{labelaff} imply that
\begin{equation} \label{labelafg}
t_k = \frac{1}{\sqrt{q_k}}, \quad k \geq 1.
\end{equation}
Therefore, the state parameters of Algorithm~\ref{labelaew} no longer contain information on the convergence guarantee but can be brought to a simpler form. It follows from \eqref{labelaeo}, \eqref{labelaep}, and \eqref{labelafg} that the auxiliary point extrapolation factor is given by
\begin{equation} \label{labelafh}
\beta_k = \frac{\sqrt{L_k + \mu_{\Psi}} - \bm{1}_{\{\bm{z}_k\}}(\bm{x}_k) \sqrt{\mu}}{\sqrt{L_{k + 1} + \mu_{\Psi}} + \sqrt{\mu}}, \quad k \geq 0.
\end{equation}
The sequence $\{t_k\}_{k \geq 0}$ does not store any relevant information and can be left out. This means that the convergence guarantee $A_k$ requires a dedicated update. A simple recursion rule follows from \eqref{labelabh}, \eqref{labelabs}, and \eqref{labeladg} as
\begin{equation} \label{labelafi}
A_{k + 1} = \frac{\sqrt{L_{k + 1} + \mu_{\Psi}}}{\sqrt{L_{k + 1} + \mu_{\Psi}} - \sqrt{\mu}} A_k \quad k \geq 0.
\end{equation}
Due to scaling invariance, we can select any pair $(A_0, \gamma_0)$ that is a positive multiple of $(1, \mu)$. For simplicity, we choose $A_0 = 1$ and $\gamma_0 = \mu$.

To reduce computational intensity, we modify subexpressions $\bm{d}_k$ and $\omega_k$ as
\begin{align}
\bm{d}_k &= \left(\sqrt{L_k + \mu_{\Psi}} - \bm{1}_{\{\bm{z}_k\}}(\bm{x}_k) \sqrt{\mu} \right) (\bm{z}_k - \bm{x}_{k - 1}) , \quad k \geq 0 , \label{labelafj} \\
\omega_k &= \frac{1}{\sqrt{L_{k + 1} + \mu_{\Psi}} + \sqrt{\mu}} , \quad k \geq 0 \label{labelafk}.
\end{align}
The local inverse condition number sequence $\{q_k\}_{k \geq 0}$ does not appear in updates \eqref{labelafi}, \eqref{labelafj}, and \eqref{labelafk}. Hence, it can also be abstracted away. The form taken by generalized ACGM in this border-case, after simplifications, is listed in Algorithm~\ref{labelafl}.

\begin{algorithm}[h!]
\caption{Border-case ACGM in extrapolated form \newline
\textbf{ACGM}($\bm{x_0}$, $L_0$, $\mu_f$, $\mu_{\Psi}$, $r_u$, $r_d$, $K$)}
\label{labelafl}
\begin{algorithmic}[1]
\STATE $\bm{x}_{-1} = \bm{x}_0, \ \bm{d}_0 = \bm{0}, \ A_0 = 1, \ \mu = \mu_f + \mu_{\Psi}$
\FOR {$k = 0, ..., K - 1$}
\STATE $\cmdaab{L} := r_d L_k$ \label{labelafm}
\LOOP
\STATE $\cmdaac{y} := \bm{x}_k + \frac{1}{\sqrt{\cmdaab{L} + \mu_{\Psi}} + \sqrt{\mu}} \bm{d}_k$ \label{labelafn}
\STATE $\cmdaac{z} :=\mbox{prox}_{\frac{1}{\cmdaab{L}} \Psi} \left(\cmdaac{y} - \frac{1}{\cmdaab{L}}\nabla f(\cmdaac{y}) \right)$ \label{labelafo}

\IF {$f(\cmdaac{z}) \leq Q_{\cmdaab{L}, \cmdaac{y}}(\cmdaac{z})$} \label{labelafp}
\STATE Break from loop \label{labelafq}
\ELSE
\STATE $\cmdaab{L} := r_u \cmdaab{L} $
\ENDIF
\ENDLOOP
\STATE $\bm{x}_{k + 1} := \argmin \{ F(\cmdaac{z}), F(\bm{x}_k) \}$ \label{labelafr}
\STATE $\bm{d}_{k + 1} := \left(\sqrt{\cmdaab{L} + \mu_{\Psi}} - \bm{1}_{\{\cmdaac{z}\}}(\bm{x}_{k + 1} ) \sqrt{\mu} \right)(\cmdaac{z} - \bm{x}_k)$
\STATE $L_{k + 1} := \cmdaab{L}$
\STATE $A_{k + 1} := \frac{\sqrt{\cmdaab{L} + \mu_{\Psi}}}{\sqrt{\cmdaab{L} + \mu_{\Psi}} - \sqrt{\mu}} A_k$
\ENDFOR
\STATE \textbf{return} $x_K$
\end{algorithmic}
\end{algorithm}

A non-monotone variant can be obtained by replacing \cmdaah{labelafl}{labelafr} with \eqref{labeladh}. The per-iteration complexity of Algorithm~\ref{labelafl} matches the one of Algorithms \ref{labeladi}~and~\ref{labelaew}.

\section{Simulation results}

\subsection{Benchmark setup}

We have tested the variants of generalized ACGM introduced in this work against the methods considered at the time of writing to be the state-of-the-art on the problem class outlined in Subsection~\ref{labelaai}. The proposed methods included in the benchmark are non-monotone ACGM (denoted as plain ACGM), monotone ACGM (MACGM), and, for strongly-convex problems, border-case non-monotone ACGM (BACGM) as well as border-case monotone ACGM (BMACGM). The state-of-the-art methods are FISTA-CP, monotone FISTA-CP (MFISTA-CP)~\cite{refaah}, AMGS~\cite{refaak}, and FISTA with backtracking line-search (FISTA-BT)~\cite{refaal}.

We have selected as test cases five synthetic instances of composite problems in the areas of statistics, inverse problems, and machine learning. Three are non-strongly convex: least absolute shrinkage and selection operator (LASSO)~\cite{refaam}, non-negative least squares (NNLS), and $l_1$-regularized logistic regression (L1LR). The other two are strongly-convex: ridge regression (RR) and elastic net (EN)~\cite{refabk}. Table~\ref{labelaab} lists the oracle functions of all above mentioned problems.

Here, the sum softplus function $\mathcal{I}(\bm{x})$, the element-wise logistic function $\mathcal{L}(\bm{x})$, and the shrinkage operator $\mathcal{T}_{\tau}(\bm{x})$ are, respectively, given for all $\bm{x} \in \mathbb{R}^n$ and $\tau > 0$ by
\begin{gather}
\mathcal{I}(\bm{x}) = \displaystyle \sum_{i = 1}^{m} \log(1 + e^{\bm{x}_i}), \\
\mathcal{L}(\bm{x})_i = \frac{1}{1 + e^{-\bm{x}_i}}, \quad i \cmdaae{m}, \\
\mathcal{T}_{\tau}(\bm{x})_j = (|\bm{x}_j| - \tau)_{+} \sgn(\bm{x}_j), \quad j \cmdaae{n}.
\end{gather}

\cmdaaj

To attain the best convergence guarantees for AMGS, Nesterov suggests in \cite{refaak} that all known global strong convexity be transferred to the simple function $\Psi$. When line-search is enabled, generalized ACGM also benefits slightly from this arrangement when \mbox{$r_u > 1$} (Theorem~\ref{labeladv}). Without line-search, the convergence guarantees of generalized ACGM do not change as a result of strong convexity transfer, in either direction. Thus, for fair comparison, we have incorporated in $\Psi$ the strongly-convex quadratic regularization term for the RR and EN problems. In the following, we describe in detail each of the five problem instances. All random variables are independent and identically distributed, unless stated otherwise.

\cmdaaf{LASSO} Real-valued matrix $\bm{A}$ is of size $m = 500$ by \mbox{$n = 500$}, with entries drawn from $\mathcal{N}(0, 1)$. Vector $\bm{b} \in \mathbb{R}^m$ has entries sampled from $\mathcal{N}(0, 9)$. Regularization parameter $\lambda_1$ is $4$. The starting point $\bm{x}_0 \in \mathbb{R}^n$ has entries drawn from $\mathcal{N}(0, 1)$.

\cmdaaf{NNLS} Sparse $m = 1000 \times n = 10000$ matrix $\bm{A}$ has approximately $10\%$ of entries, at random locations, non-zero. The non-zero entries are drawn from $\mathcal{N}(0, 1)$ after which each column $j \cmdaae{n}$ is scaled independently to have an $l_2$ norm of $1$. Starting point $\bm{x}_0$ has $10$ entries at random locations all equal to $4$ and the remainder zero. Vector $\bm{b}$ is obtained from $\bm{b} = \bm{A} \bm{x}_0 + \bm{z}$, where $\bm{z}$ is standard Gaussian noise.

\cmdaaf{L1LR} Matrix $\bm{A}$ has $m = 200 \times n = 1000$ entries sampled from $\mathcal{N}(0, 1)$, $\bm{x}_0$ has exactly $10$ non-zero entries at random locations, each entry value drawn from $\mathcal{N}(0, 225)$, and $\lambda_1 = 5$. Labels $\bm{y}_i \in \{0, 1\}, i \cmdaae{m}$ are selected with probability $\mathbb{P}(\bm{Y}_i = 1) = \mathcal{L}(\bm{A} \bm {x})_i$.

\cmdaaf{RR} Dimensions are $m = 500 \times n = 500$. The entries of matrix $\bm{A}$, vector $\bm{b}$, and starting point $\bm{x}_0$ are drawn from $\mathcal{N}(0, 1)$, $\mathcal{N}(0, 25)$, and $\mathcal{N}(0, 1)$, respectively. Regularizer $\lambda_2$ is given by $10^{-3} (\sigma_{max}(\bm{A}))^2$, where $\sigma_{max}(\bm{A})$ is the largest singular value of $\bm{A}$.

\cmdaaf{EN} Matrix $\bm{A}$ has $m = 1000 \times n = 500$ entries sampled from $\mathcal{N}(0, 1)$. Starting point $\bm{x}_0$ has $20$ non-zero entries at random locations, each entry value drawn from $\mathcal{N}(0, 1)$. Regularization parameter $\lambda_1$ is obtained according to \cite{refabl} as $\lambda_1 = 1.5 \sqrt{2 \log(n)}$ and $\lambda_2$ is the same as in RR, $\lambda_2 = 10^{-3} (\sigma_{max}(\bm{A}))^2$.

The Lipschitz constant $L_f$ is given by $(\sigma_{max}(\bm{A}))^2$ for all problems except for L1LR where it is $\frac{1}{4}(\sigma_{max}(\bm{A}))^2$. Strongly convex problems RR and EN have $\mu = \mu_{\Psi} = \lambda_2$ and inverse condition number $q = \mu / (L_f + \mu_{\Psi}) = 1 / 1001$.

To be able to benchmark against FISTA-CP and FISTA-BT, which lack fully adaptive line-search, we have set $L_0 = L_f$ for all tested algorithms, thus giving FISTA-CP and FISTA-BT an advantage over the proposed methods. To highlight the differences between ACGM and BACGM, we ran ACGM and MACGM with parameters $A_0 = 0$ and $\gamma_0 = 1$.

Despite the problems differing in structure, the oracle functions have the same computational costs. We consider one matrix-vector multiplication to cost 1 WTU. Consequently, for all problems, we have $t_f = 1$~WTU, $t_g = 2$~WTU, and $t_{\Psi} = t_p = 0$~WTU.

The line-search parameters were selected according to the recommendation given in \cite{refaaf}. For AMGS and FISTA-BT we have $r_u\cmdaad{AMGS} = r_u\cmdaad{FISTA} = 2.0$ and $r_d\cmdaad{AMGS} = 0.9$. The variants of generalized ACGM and AMGS are the only methods included in the benchmark that are equipped with fully adaptive line-search. We have decided to select $r_d\cmdaad{ACGM}$ to ensure that ACGM and AMGS have the same overhead. We formally define the line-search overhead of method $\mathcal{M}$, denoted by $\Omega\cmdaad{\mathcal{M}}$, as the average computational cost attributable to backtracks per WTU of advancement. Assuming that the LCEs hover around a fixed value (Subsection~\ref{labeladz}), we thus have that
\begin{align}
\Omega\cmdaad{AMGS} &= -\frac{( 2 t_g + t_p) \log(r_d\cmdaad{AMGS})}{2 (t_g + t_p) \log(r_u\cmdaad{AMGS})}, \label{labelafs}\\
\Omega\cmdaad{ACGM} &= -\frac{(t_f + t_g + t_p) \log(r_d\cmdaad{ACGM})}{ (t_g + t_p) \log(r_u\cmdaad{ACGM})}. \label{labelaft}
\end{align}
From \eqref{labelafs} and \eqref{labelaft} we have that $r_d\cmdaad{ACGM} = (r_d\cmdaad{AMGS})^{\frac{2}{3}}$, with no difference for border-case or monotone variants.

For measuring ISDs, we have computed beforehand an optimal point estimate $\hat{\bm{x}}^*$ for each problem instance. Each $\hat{\bm{x}}^*$ was obtained as the main iterate after running MACGM for $5000$ iterations with parameters $A_0 = 0$, $\gamma_0 = 1$, $L_0 = L_f$, and aggressive search parameters $r_d = 0.9$ and $r_u = 2.0$.

\subsection{Non-strongly convex problems}

The convergence results for LASSO, NNLS, and L1LR are shown in Figure~\ref{labelaac}. The LCE variation during the first 200 WTU is shown in Figure~\ref{labelaad}. For NNLS, floating point precision was exhausted after 100 WTU and the LCE variation was only plotted to this point (Figure~\ref{labelaad}(b)). In addition, the average LCEs are listed in Table~\ref{labelaae}.

Both variants of ACGM outperform in iterations and especially in WTU the competing methods in each of these problem instances. Even though for LASSO and NNLS, the iteration convergence rate of AMGS is slightly better in the beginning (Figures \ref{labelaac}(a) and \ref{labelaac}(c)), AMGS lags behind afterwards and, when measured in terms of computational convergence rate, has the poorest performance among the methods tested (Figures \ref{labelaac}(b), \ref{labelaac}(d), and \ref{labelaac}(f)). FISTA-BT produces the same iterates as FISTA-CP, as theoretically guaranteed in the non-strongly convex case for $L_0 = L_f$.

The overall superiority of ACGM and MACGM can be attributed to the effectiveness of line-search. Interestingly, ACGM manages to surpass FISTA-CP and MFISTA-CP even when the latter are supplied with the exact value of the global Lipschitz constant. This is because ACGM is able to accurately estimate the local curvature, which is often well below $L_f$. For the L1LR problem, where the smooth part $f$ is not the square of a linear function, the local curvature is substantially lower than the global Lipschitz constant with LCEs hovering around one fifth of $L_f$ (Figure~\ref{labelaad}(c)). One would expect AMGS to be able to estimate local curvature as accurately as ACGM. This is does not happen due to AMGS's reliance on a ``damped relaxation condition''~\cite{refaak} line-search stopping criterion. For LASSO and NNLS, the average LCE of AMGS is actually above $L_f$. ACGM has an average LCE that is roughly two thirds that of AMGS on these problems whereas for L1LR the average is more than three times lower than AMGS. The difference between the LCE averages of ACGM and MACGM is negligible.

Indeed, monotonicity, as predicted, does not alter the overall iteration convergence rate and has a stabilizing effect. MACGM overshoots do have a negative but limited impact on the computational convergence rate. We have noticed in our simulations that overshoots occur less often for larger problems, such as the tested instance of NNLS.

\subsection{Strongly convex problems}

The convergence results for RR and EN are shown in Figures \ref{labelaaf}(a), \ref{labelaaf}(b), \ref{labelaaf}(c), and \ref{labelaaf}(d). The LCE variation is shown in Figure~\ref{labelaag} with LCE averages listed in Table~\ref{labelaah}.

Strongly convex problems have a unique optimum point and accelerated first-order schemes are guaranteed to find an accurate estimate of it in domain space (see \cite{refaae} for a detailed analysis). Along with Theorem~\ref{labeladv}, it follows that
\begin{equation}
A_0 (F(\bm{x}_0) - F(\hat{\bm{x}}^*)) + \frac{\gamma_0}{2}\| \bm{x}_0 - \hat{\bm{x}}^* \|_2^2 \simeq \Delta_0.
\end{equation}
Thus, we can display accurate estimates of ISDUBs in \eqref{labelaam}, of the form $U_k \operatorname{\overset{def}{=}} \Delta_0 / A_k$, $k \geq 1$, for methods that maintain convergence guarantees at runtime. These are shown in Figures \ref{labelaaf}(e) and \ref{labelaaf}(f) as upper bounds indexed in WTU.
\clearpage
\cmdaak
\clearpage
\cmdaal

For the smooth RR problem, the effectiveness of each tested algorithm is roughly given by the increase rate of the accumulated weights (Figures \ref{labelaaf}(a) and \ref{labelaaf}(c)). In iterations, AMGS converges the fastest. However, in terms of WTU usage, it is the least effective of the methods designed to deal with strongly convex objectives. The reasons are the high cost of its iterations, its low asymptotic rate compared to ACGM and FISTA-CP, and the stringency of its damped relaxation criterion that results in higher LCEs (on average) than ACGM (Figure~\ref{labelaag}(a) and Table~\ref{labelaah}). The computational convergence rate of BACGM is the best, followed by ACGM and FISTA-CP. This does not, however, correspond to the upper bounds (Figure \ref{labelaaf}(e)). While BACGM produces the largest accumulated weights $A_k$, the high value of the ISD term in $\Delta_0$ causes BACGM to have poorer upper bounds than ACGM, except for the first iterations. In fact, the effectiveness of BACGM on this problem is exceptional, partly due to the regularity of the composite gradients. This regularity also ensures the monotonicity of BACGM, ACGM, and FISTA-CP. FISTA-BT does not exhibit linear convergence on this problem and it is even slower than AMGS after 500 WTU despite its lower line-search overhead and advantageous parameter choice $L_0 = L_f$.
\clearpage
\cmdaam
\clearpage
\cmdaan

On the less regular EN problem, ACGM leads all other methods in terms of both iteration and computational convergence rates (Figures \ref{labelaaf}(b) and \ref{labelaaf}(d)). The advantage of ACGM, especially over BACGM, is accurately reflected in the upper bounds (Figure~\ref{labelaaf}(f)). However, convergence is much faster than the upper bounds would imply. Even FISTA-BT has a competitive rate, due to the small number of iterations ($150$) needed for high accuracy results. The ineffectiveness of AMGS on this problem is mostly due to its high LCEs (Figure~\ref{labelaag}(b) and Table~\ref{labelaah}). The proposed ACGM and variants show comparable average LCEs. Here as well, monotonicity has a stabilizing effect and does not have a significant impact on the computational convergence rate.

\section{Discussion} \label{labelafu}

Our simulation results suggest that enforcing monotonicity in ACGM is generally beneficial in large-scale applications. It leads to a more predictable convergence rate and, provided that the number of overshoots per iteration is small, monotonicity has a negligible impact on the computational convergence rate as well. Our experimental results indicate that the frequency of overshoots generally decreases with problem size.

From a theoretical standpoint, the proposed method can be viewed as a unification of FGM and FISTA, in their most common forms. Specifically, the fixed-step variant (\mbox{$L_k = L_f$} for all $k \geq 0$) of ACGM in extrapolated form (Algorithm~\ref{labelaew}) is equivalent to both the monotone and non-monotone variants of FISTA-CP with the theoretically optimal step size $\tau\cmdaad{FISTA-CP} = 1 / L_f$. Moreover, when $\mu = 0$, the original non-monotone fixed-step ACGM coincides with the original formulation of FISTA in~\cite{refaal}. Adding monotonicity yields MFISTA~\cite{refaba}. Also for $\mu = 0$, but without the fixed-step restriction, the original non-monotone ACGM in estimate sequence form reduces to the robust FISTA-like algorithm in~\cite{refaay}, whereas in extrapolated form it is a valuable simplification of the method introduced in~\cite{refaax}.
The original backtracking FISTA can be obtained in the same way as the variants in \cite{refaay, refaax} by further removing the ratio $L_{k + 1} / L_k$ from the $t_{k + 1}$ update in line~\ref{labelaez} of Algorithm~\ref{labelaew} (noting that $\mu_{\Psi} = 0$ for $\mu = 0$). Since backtracking FISTA relies on the assumption that this ratio is never less that $1$, the convergence analysis remains unaffected.

When dealing with differentiable objectives, we can assume without loss of generality that $\Psi(\bm{x}) = 0$ for all $\bm{x} \in \mathbb{R}^n$. In what follows, we consider generalized non-monotone fixed-step ACGM in estimate sequence form, unless stated otherwise. By substituting the local upper bound functions $u_{k + 1}$ at every iteration $k \geq 0$ with any functions that produce iterates satisfying the descent condition, which means in this context that
\begin{equation} \nonumber
f(\bm{x}_{k + 1}) \leq f(\bm{y}_{k + 1}) - \frac{1}{2 L_{k + 1}} \| \nabla f(\bm{y}_{k + 1}) \|_2^2,
\end{equation}
where $\bm{x}_{k + 1}$ is given by \cmdaai{labelacc}{labelaci}, we obtain the ``general scheme of optimal method'' in \cite{refaae}. Both the monotone and non-monotone variants adhere to this scheme. The correspondence between Nesterov's notation in \cite{refaae} and ours is for all $k \geq 0$ given by:
\begin{equation} \nonumber
\begin{gathered}
\lambda_k\cmdaad{FGM} = \frac{A_0\cmdaad{ACGM}}{A_k\cmdaad{ACGM}}, \quad
\phi_k\cmdaad{FGM}(\bm{x}) = \frac{1}{A_k\cmdaad{ACGM}} \psi_k\cmdaad{ACGM}(\bm{x}), \quad \bm{x} \in \mathbb{R}^n, \\
\bm{y}_k\cmdaad{FGM} = \bm{y}_{k + 1}\cmdaad{ACGM}, \quad
\alpha_k\cmdaad{FGM} = \frac{a_{k + 1}\cmdaad{ACGM}}{A_{k + 1}\cmdaad{ACGM}} = \frac{1}{t_{k + 1}\cmdaad{ACGM}}, \quad
\gamma_k\cmdaad{FGM} = \frac{\gamma_k\cmdaad{ACGM}}{A_k\cmdaad{ACGM}}.
\end{gathered}
\end{equation}

The remaining state parameters are identical. Note that FGM makes the assumption that $A_0\cmdaad{ACGM} > 0$, which is incompatible with the original specification of ACGM in \cite{refaaz}. With the above assumption, the non-monotone variant (Algorithm~\ref{labeladi}) is in fact identical to the ``constant step scheme I'' in \cite{refaae}. Similarly, the extrapolated form of fixed-step non-monotone ACGM (Algorithm~\ref{labelaew}) corresponds exactly to the ``constant step scheme II'' in \cite{refaae} while fixed-step border-case non-monotone ACGM (Algorithm~\ref{labelafl}) is identical to the ``constant step scheme III'' in \cite{refaae}. We note that for both the RR and EN problems, regardless of the actual performance of BACGM, the convergence guarantees of BACGM are poorer than those of ACGM with $A_0 = 0$. This discrepancy in guarantees is supported theoretically because, in the most common applications, the ISD term in $\Delta_0$ is large compared to the DST. This extends to the fixed-step setup and challenges the notion found in the literature (e.g., \cite{refabm}) that for strongly-convex functions, FGM and FISTA-CP are momentum methods that take the form of the ``constant step scheme III'' in \cite{refaae}. In fact, the border-case form may constitute the poorest choice of parameters $A_0$ and $\gamma_0$ in many applications. Indeed, the worst-case results in Theorem~\ref{labeladv} favor $A_0 = 0$.

The FGM variant in~\cite{refaaw} is a particular case of non-monotone ACGM with variable step size (Algorithm~\ref{labeladi}) when the objective is non-strongly convex ($\mu = 0$) and the step size search parameters are set to $r_u\cmdaad{ACGM} = 2$ and $r_d\cmdaad{ACGM} = 0.5$. The notation correspondence is as follows:
\begin{equation} \nonumber
\bm{x}_{k + 1, i}\cmdaad{FGM} = \bm{x}_{k + 1}\cmdaad{ACGM}, \quad \bm{y}_{k, i}\cmdaad{FGM} = \bm{y}_{k + 1}\cmdaad{ACGM}, \quad a_{k, i}\cmdaad{FGM} = \cmdaab{a}\cmdaad{ACGM}, \quad 2^{i} L_f = \cmdaab{L}\cmdaad{ACGM}.
\end{equation}
The remaining parameters are identical.

Thus, by relaxing the assumption that $A_0 = 0$, we have devised a generalized variant of ACGM that effectively encompasses FGM~\cite{refaae}, with its recently introduced variant~\cite{refaaw}, as well as the original FISTA~\cite{refaal}, including its adaptive step-size variants~\cite{refaay,refaax}, the monotone version MFISTA~\cite{refaba}, and the strongly convex extensions FISTA-CP~\cite{refaah} and scAPG~\cite{refaap}. A summary of how the above first-order methods relate to generalized ACGM is given in Table~\ref{labelafv}.

\begin{table*}[t!]
\centering \small
\caption{FGM and FISTA, along with their common variants, can be considered instances of generalized ACGM with certain restrictions applied.}
\label{labelafv}
\begin{tabular}{@{}m{21mm} m{14mm} m{10mm} m{10mm} m{10mm} m{10mm} m{12mm} m{14mm}@{}} \toprule
Algorithm & \multicolumn{7}{c}{Restriction} \\ & Smooth objective & $\mu = 0$ & $\mu > 0$ & $A_0 = 0$ & $A_0 > 0$ & Fixed step size & Non-monotone \\ \midrule
FGM~\cite{refaae} & yes & no & no & no & yes & yes & yes \\
FGM~\cite{refaaw} & yes & yes & no & unclear & unclear & no & yes \\
FISTA~\cite{refaal} & no & yes & no & yes & no & partial & yes \\
FISTA~\cite{refaax,refaay} & no & yes & no & yes & no & no & yes \\
MFISTA~\cite{refaba} & no & yes & no & yes & no & yes & no \\
scAPG~\cite{refaap} & no & no & yes & no & yes & no & yes \\
FISTA-CP~\cite{refaah} & no & no & no & no & no & yes & no \\
\bottomrule
\end{tabular}
\end{table*}

Due to its adaptivity, generalized ACGM is not limited to the composite problem framework in Subsection~\ref{labelaai}. It is also guaranteed to converge on problems where the gradient of $f$ is \emph{not globally} Lipschitz continuous. Constituent function $f$ needs to have Lipschitz gradient only in the area explored by the algorithm.

In terms of usability, generalized ACGM does not require a priori knowledge of Lipschitz constant $L_f$, or a lower estimate of it, beforehand. Thus, the proposed method can be utilized \emph{without any quantitative knowledge of the problem}. Lack of information does not hinder generalized ACGM more than any other primal first-order method while additional information, such as the values of strong convexity parameters $\mu_f$ and $\mu_{\Psi}$ or even an accurate estimate $L_0$ of the curvature around $\bm{x}_0$, leads to a state-of-the-art performance increase unsurpassed for its class.

\bibliographystyle{IEEEtran}

\begin{thebibliography}{10}

\bibitem{refaab}
A.~Nemirovski and D.-B. Yudin,
\emph{Problem complexity and method efficiency in optimization}.
John Wiley \& Sons, New York, NY, USA, 1983.

\bibitem{refaac}
Y.~Nesterov,
``Subgradient methods for huge-scale optimization problems,''
\emph{Math. Program., Ser. A}, vol. 146, no. 1-2, pp. 275--297, 2014.

\bibitem{refaad}
------,
``A method of solving a convex programming problem with convergence rate $\mathcal{O}(1/k^2)$,''
\emph{Dokl. Math.}, vol.~27, no.~2, pp. 372--376, 1983.

\bibitem{refaae}
------, \emph{Introductory Lectures on Convex Optimization. Applied Optimization, vol. 87}.
Kluwer Academic Publishers, Boston, MA, USA, 2004.

\bibitem{refaaf}
S.~R. Becker, E.~J. Cand{\`e}s, and M.~C. Grant,
``Templates for convex cone problems with applications to sparse signal recovery,''
\emph{Math. Program. Comput.}, vol.~3, no.~3, pp. 165--218, 2011.

\bibitem{refaag}
S.~Bubeck \emph{et~al.},
``Convex optimization: Algorithms and complexity,''
\emph{Found. Trends Mach. Learn.}, vol.~8, no. 3-4, pp. 231--357, 2015.

\bibitem{refaah}
A.~Chambolle and T.~Pock,
``An introduction to continuous optimization for imaging,''
\emph{Acta Numer.}, vol.~25, pp. 161--319, 2016.

\bibitem{refaai}
N.~Parikh, S.~P. Boyd \emph{et~al.},
``Proximal algorithms,''
\emph{Found. Trends Optim.}, vol.~1, no.~3, pp. 127--239, 2014.

\bibitem{refaaj}
K.~Slavakis, G.~B. Giannakis, and G.~Mateos,
``Modeling and optimization for big data analytics:(statistical) learning tools for our era of data deluge,''
\emph{{IEEE} Signal Process. Mag.}, vol.~31, no.~5, pp. 18--31, Sep. 2014.

\bibitem{refaak}
Y.~Nesterov,
``Gradient methods for minimizing composite functions,''
\emph{Universit\'e catholique de Louvain, CORE Discussion Papers}, 2007/76, Sep. 2007.

\bibitem{refaal}
A.~Beck and M.~Teboulle,
``A fast iterative shrinkage-thresholding algorithm for linear inverse problems,''
\emph{{SIAM} J. Imaging Sci.}, vol.~2, no.~1, pp. 183--202, 2009.

\bibitem{refaam}
R.~Tibshirani,
``Regression shrinkage and selection via the lasso,''
\emph{J. R. Stat. Soc. Ser. B. Methodol.}, vol.~58, no.~1, pp. 267--288, 1996.

\bibitem{refaan}
J.~Mairal,
``Optimization with first-order surrogate functions,''
in \emph{ICML}, 2014, Atlanta, Georgia, USA, pp. 73--81.

\bibitem{refaao}
H.~Lin, J.~Mairal, and Z.~Harchaoui,
``A universal catalyst for first-order optimization,''
in \emph{NIPS}, Dec. 2015, Montreal, Canada, pp. 3384--3392.

\bibitem{refaap}
Q.~Lin and L.~Xiao,
``An adaptive accelerated proximal gradient method and its homotopy continuation for sparse optimization,''
in \emph{ICML}, 2014, pp. 73--81.

\bibitem{refaaq}
B.~O'Donoghue and E.~Cand{\`e}s,
``Adaptive restart for accelerated gradient schemes,''
\emph{Found. Comput. Math.}, vol.~15, no.~3, pp. 715--732, 2015.

\bibitem{refaar}
M.~Yamagishi and I.~Yamada,
``Over-relaxation of the fast iterative shrinkage-thresholding algorithm with variable stepsize,''
\emph{Inverse Problems}, vol.~27, no.~10, p. 105008, Sep. 2011.

\bibitem{refaas}
W.~W. Hager, D.~T. Phan, and H.~Zhang, ``Gradient-based methods for sparse recovery,''
\emph{{SIAM} J. Imaging Sci.}, vol.~4, no.~1, pp. 146--165, 2011.

\bibitem{refaat}
Z.~Wen, W.~Yin, D.~Goldfarb, and Y.~Zhang,
``A fast algorithm for sparse reconstruction based on shrinkage, subspace optimization, and continuation,''
\emph{{SIAM} J. Sci. Comput.}, vol.~32, no.~4, pp. 1832--1857, 2010.

\bibitem{refaau}
Z.~Wen, W.~Yin, H.~Zhang, and D.~Goldfarb,
``On the convergence of an active-set method for $l_1$ minimization,''
\emph{Optim. Methods Software}, vol.~27, no.~6, pp. 1127--1146, 2012.

\bibitem{refaav}
S.~J. Wright, R.~D. Nowak, and M.~A. Figueiredo,
``Sparse reconstruction by separable approximation,''
\emph{{IEEE} Trans. Signal Process.}, vol.~57, no.~7, pp. 2479--2493, 2009.

\bibitem{refaaw}
Y.~Nesterov and S.~Stich,
``Efficiency of accelerated coordinate descent method on structured optimization problems,''
\emph{Universit\'e catholique de Louvain, CORE Discussion Papers}, 2016/03, Feb. 2016.

\bibitem{refaax}
K.~Scheinberg, D.~Goldfarb, and X.~Bai,
``Fast first-order methods for composite convex optimization with backtracking,''
\emph{Found. Comput. Math.}, vol.~14, no.~3, pp. 389--417, 2014.

\bibitem{refaay}
M.~I. Florea and S.~A. Vorobyov,
``A robust {FISTA}-like algorithm,''
in \emph{Proc. of {IEEE} Intern. Conf. on Acoustics, Speech and Signal Processing (ICASSP)}, Mar. 2017, New Orleans, USA, pp. 4521--4525.

\bibitem{refaaz}
------,
``An accelerated composite gradient method for large-scale composite objective problems,''
\emph{{IEEE} Trans. Signal Process.}, vol.~67, no.~2, pp. 444--459, Jan. 2019.

\bibitem{refaba}
A.~Beck and M.~Teboulle,
``Fast gradient-based algorithms for constrained total variation image denoising and deblurring problems,''
\emph{{IEEE} Trans. Image Process.}, vol.~18, no.~11, pp. 2419--2434, 2009.

\bibitem{refabb}
A.~Chambolle and C.~Dossal,
``On the convergence of the iterates of the Fast Iterative Shrinkage/Thresholding Algorithm,''
\emph{J. Optim. Theory Appl}, vol. 166, no.~3, pp. 968--982, 2015.

\bibitem{refabc}
B.~Polyak,
\emph{Introduction to optimization (translated from Russian)}.
Translations Series in Mathematics and Engineering, Optimization Software, New York, NY, USA, 1987.

\bibitem{refabd}
A.~Brown and M.~C. Bartholomew-Biggs,
``Some effective methods for unconstrained optimization based on the solution of systems of ordinary differential equations,''
\emph{J. Optim. Theory Appl}, vol.~62, no.~2, pp. 211--224, 1989.

\bibitem{refabe}
S.~Boyd, N.~Parikh, E.~Chu, B.~Peleato, and J.~Eckstein,
``Distributed optimization and statistical learning via the alternating direction method of multipliers,''
\emph{Found. Trends Mach. Learn.}, vol.~3, no.~1, pp. 1--122, 2011.

\bibitem{refabf}
J.~M. Ortega and W.~C. Rheinboldt,
\emph{Iterative solution of nonlinear equations in several variables}.
Society for Industrial and Applied Mathematics, Philadelphia, PA, USA, 1970.

\bibitem{refabg}
K.~Lange,
\emph{MM optimization algorithms}.
Society for Industrial and Applied Mathematics, Philadelphia, PA, USA, 2016.

\bibitem{refabh}
R.~T. Rockafellar,
\emph{Convex Analysis}.
Princeton University Press, Princeton, NJ, USA, 1970.

\bibitem{refabi}
M.~I. Florea,
``Constructing accelerated algorithms for large-scale optimization,''
Ph.D. dissertation, Aalto University, School of Electrical Engineering, Helsinki, Finland, Oct. 2018.

\bibitem{refabj}
J.~L. Hennessy and D.~A. Patterson,
\emph{Computer Architecture: A Quantitative Approach}, 5th~ed.
Morgan Kaufmann Publishers, San Francisco, CA, USA, 2011.

\bibitem{refabk}
H.~Zou and T.~Hastie,
``Regularization and variable selection via the elastic net,''
\emph{J. R. Stat. Soc. Ser. B. Methodol.}, vol.~67, no.~2, pp. 301--320, 2005.

\bibitem{refabl}
T.~Hastie, R.~Tibshirani, and M.~Wainwright,
\emph{Statistical Learning with Sparsity: The Lasso and Generalizations}.
CRC Press, 2015.

\bibitem{refabm}
W.~Su, S.~Boyd, and E.~J. Cand{\`e}s,
``A differential equation for modeling {N}esterov's accelerated gradient method: Theory and insights,''
\emph{J. Mach. Learn. Res.}, vol.~17, pp. 1--43, 2016.
\end{thebibliography}

\end{document}